\documentclass[11pt,a4paper]{amsart}

 \usepackage{amsfonts,amsmath,amscd,amssymb,amsbsy,amsthm,amstext,amsopn}
 \usepackage{fullpage,mathrsfs,subfigure}
 \usepackage[all,dvips,arc,curve,color,frame]{xy}
 \usepackage{pstricks,pst-node}

 \numberwithin{equation}{section}
\newtheorem{theorem}{Theorem}[section]
\newtheorem{proposition}[theorem]{Proposition}
\newtheorem{lemma}[theorem]{Lemma}
\newtheorem{corollary}[theorem]{Corollary}

\theoremstyle{definition}

\newtheorem{example}[theorem]{Example}

\theoremstyle{remark}
\newtheorem{remark}[theorem]{Remark}

\newcommand{\kk}{\ensuremath{\Bbbk}} 

\newcommand{\NN}{\ensuremath{\mathbb{N}}} 
 
\newcommand{\QQ}{\ensuremath{\mathbb{Q}}} 
 
\newcommand{\ZZ}{\ensuremath{\mathbb{Z}}} 
\newcommand{\one}{\ensuremath{(\mathrm{i})}}
\newcommand{\two}{\ensuremath{(\mathrm{ii})}}
\newcommand{\three}{\ensuremath{(\mathrm{iii})}}
\newcommand{\four}{\ensuremath{(\mathrm{iv})}}

\newcommand{\head}{\operatorname{h}}

\newcommand{\inc}{\operatorname{inc}}

\newcommand{\mon}{\operatorname{mon}}

\newcommand{\pic}{\operatorname{pic}}

\newcommand{\supp}{\operatorname{supp}}
\newcommand{\tail}{\operatorname{t}}

\newcommand{\Coh}{\operatorname{Coh}}

\newcommand{\End}{\operatorname{End}}

\newcommand{\GL}{\operatorname{GL}}

\newcommand{\Hom}{\operatorname{Hom}} 
 
\newcommand{\Irr}{\operatorname{Irr}} 
\newcommand{\Ker}{\operatorname{Ker}}

\newcommand{\Pic}{\operatorname{Pic}}
\newcommand{\Proj}{\operatorname{Proj}}
 
\newcommand{\SL}{\operatorname{SL}} 
\newcommand{\Spec}{\operatorname{Spec}}

\newcommand{\Wt}{\operatorname{Wt}}

 \newcommand{\modA}{\operatorname{mod}(\ensuremath{A})}

\newcommand{\git}{\ensuremath{/\!\!/\!}}
\renewcommand{\div}{\operatorname{div}} 

 \newcommand{\Irrsp}{\Irr^\mathrm{sp}}
 \newcommand{\ghilb}{\ensuremath{G}\operatorname{-Hilb}}
 \newcommand{\Gammahilb}{\ensuremath{\Gamma}\operatorname{-Hilb}}
 \newcommand{\Per}{\ensuremath{^{-1}}\operatorname{Per}\bigl(X/(\mathbb{A}^2_\kk/G)\bigr)}

\begin{document}
 \bibliographystyle{plain} 
 \title[The Special McKay correspondence]{The Special McKay
   correspondence as an equivalence of derived categories}
 \author{Alastair Craw}
 \address{Department of Mathematics, Glasgow  University, Glasgow G12 8QW, Scotland.}
 \email{craw@maths.gla.ac.uk}
 
%\thanks{The author was partially supported by EPSRC}

  \begin{abstract}
    We give a new moduli construction of the minimal resolution of  
    the singularity of type $\frac{1}{r}(1,a)$ by introducing the
    Special McKay quiver.  To demonstrate that our construction trumps
    that of the $G$-Hilbert scheme, we show that the induced
    tautological line bundles freely generate the bounded derived
    category of coherent sheaves on $X$ by establishing a suitable
    derived equivalence. This gives a moduli construction of the
    Special McKay correspondence for abelian subgroups of $\GL(2)$. 
 \end{abstract}

 \maketitle

 \section{Introduction}
 For a finite subgroup $G\subset \SL(2,\kk)$, the McKay correspondence
 establishes an equivalence between the geometry of the minimal
 resolution $X$ of $\mathbb{A}^2_\kk/G$ and the $G$-equivariant
 geometry of $\mathbb{A}^2_\kk$.  More precisely, following the
 description by Ito--Nakamura~\cite{ItoNakamura} of $X$ as the
 $G$-Hilbert scheme, Kapranov--Vasserot~\cite{KapranovVasserot} used
 the resulting universal family to establish an equivalence between
 the bounded derived category of coherent sheaves on $X$ and the
 bounded derived category of finitely generated modules over the
 skew group algebra $\kk[x,y]*G$. For a finite subgroup $G\subset \GL(2,\kk)$ that is not
 special-linear, however, the $G$-Hilbert scheme has too many
 tautological bundles, so this moduli description of the McKay
 correspondence cannot hold without some redundancy.  Nevertheless,
 Wunram~\cite{Wunram} constructed an integral basis of the
 Grothendieck group of vector bundles on $X$ that is indexed by the
 trivial representation and the so-called \emph{special}
 representations of $G$.  Ishii~\cite{Ishii} subsequently employed the
 universal family for the $G$-Hilbert scheme to establish a fully
 faithful functor from the bounded derived category of compactly
 supported coherent sheaves on $X$ to the bounded derived category of
 finitely generated nilpotent modules over $\kk[x,y]*G$. 
 
 This article adopts a new approach for a finite abelian subgroup
 $G\subset \GL(2,\kk)$ by introducing a new moduli construction of the minimal resolution $X$ of $\mathbb{A}^2_\kk/G$.
 To begin, we fix a collection of line bundles
 $L_0=\mathscr{O}_X,L_1,\dots,L_\ell$ on $X$ that form an integral basis of the Grothendieck
 group of vector bundles.  The algebra
 $\End\bigl( \bigoplus_{0\leq i\leq \ell} L_i \bigr)$
 is isomorphic to the quotient $\kk Q/R$ of
 the path algebra of a quiver $Q$ modulo an ideal of relations $R$,
 where the pair $(Q,R)$ is the bound quiver of sections of the
 collection $\underline{\mathscr{L}}=(L_0,\dots,L_\ell)$ as defined by
 Craw--Smith~\cite{CrawSmith}. Since $(Q,R)$ can be
 characterised in terms of Wunram's special representations, we call $(Q,R)$ the \emph{bound Special McKay quiver}.  The main result of this article is the following.
 
 \begin{theorem} 
 \label{thm:main}
 Let $G\subset \GL(2,\kk)$ be a finite abelian subgroup with bound
 Special McKay quiver $(Q,R)$.  The minimal resolution $X$ of
 $\mathbb{A}^2_\kk/G$ is isomorphic to the fine moduli space
 $\mathcal{M}_{\vartheta}$ of $\vartheta$-stable representations
 of $(Q,R)$ for a given $\vartheta$, with tautological bundle
 $\bigoplus_{0\leq i\leq \ell} L_i$.
\end{theorem}

The proof of this result has two parts. The first, geometric part extends the construction of Craw--Smith~\cite{CrawSmith} by defining the morphism $\varphi_{\vert
\underline{\mathscr{L}}\vert}\colon X \longrightarrow \vert \underline{\mathscr{L}}\vert$ to
the multigraded linear series associated to  the sequence of line bundles $\underline{\mathscr{L}}=
(L_0,\dots,L_\ell)$ on $X$. The toric variety $\vert \underline{\mathscr{L}}\vert$ is defined to be the fine moduli space of $\vartheta$-stable representations of the quiver of sections $Q$ of $\underline{\mathscr{L}}$ for a given dimension vector and stability condition $\vartheta$, and it contains $\mathcal{M}_{\vartheta}(Q,R)$ as a closed subscheme. We prove that $\varphi_{\vert
\underline{\mathscr{L}}\vert}$ is a closed immersion and, moreover, we identify explicitly the image as a toric subvariety $\mathbb{V}(I_Q)\git_\vartheta T_Q$ of $\mathcal{M}_{\vartheta}(Q,R)$.

The second, algebraic part establishes a link between the bound quiver of sections $(Q,R)$ of $\underline{\mathscr{L}}$ and the bound McKay quiver of $G\subset
\GL(2,\kk)$. If the given two-dimensional representation of $G$
decomposes into irreducible representations as
$\mathbb{A}^2_\kk=\rho_1\oplus \rho_2$, then the McKay quiver
$Q^\prime$ is the quiver whose vertices are indexed by the irreducible
representations, and where each vertex $\rho$ admits two incoming
arrows: one arrow from vertex $\rho\otimes \rho_1$; and the second
from vertex $\rho\otimes\rho_2$. The ideal of relations
$R^\prime\subset \kk Q^\prime$ ensures that $\kk
Q^\prime/R^\prime$ is isomorphic to the skew group algebra $\kk[x,y]*G$, and the pair $(Q^\prime,
R^\prime)$ is the bound McKay quiver of $G\subset \GL(2,\kk)$.  We
show that $(Q^\prime,R^\prime)$ coincides with the bound quiver of
sections of the tautological line bundles on the $G$-Hilbert scheme
$\ghilb(\mathbb{A}^2_\kk)$ and, moreover, that the isomorphism between $\ghilb(\mathbb{A}^2_\kk)$ and the minimal resolution $X$ identifies the tautological bundles
on $\ghilb(\mathbb{A}^2_\kk)$ indexed by the trivial and special
representations of $G$ with the integral basis $L_0,L_1,\dots, L_\ell$
of the Grothendieck group of $X$. This implies in particular that $\End\bigl( \bigoplus_{0\leq i\leq \ell} L_i \bigr)$
is isomorphic to a subalgebra of $\kk
Q^\prime/R^\prime$. This description, coupled with our understanding of the relations $R^\prime$ in the McKay quiver, provides just enough information about $R$ to show that the moduli space $\mathcal{M}_{\vartheta}(Q,R)$ coincides with the image $\mathbb{V}(I_Q)\git_\vartheta T_Q$ of the morphism $\varphi_{\vert
\underline{\mathscr{L}}\vert}$.

 Theorem~\ref{thm:main} provides a moduli description for the following derived
 category version of the \emph{Special McKay correspondence} that is essentially due to Van den Bergh~\cite{VDB}. 
 
 \begin{theorem}
 \label{thm:mainderived}
 Let $\mathscr{V}:= \bigoplus_{0\leq i\leq \ell} L_i$ denote the
 tautological vector bundle on $\mathcal{M}_\vartheta(Q,R)$. Then
 \[
\mathbf{R}\!\Hom(\mathscr{V}, -) \colon
 D^b\big(\!\Coh(X)\big)\rightarrow D^b\big(\!\modA\big)
 \]
 is a derived equivalence between the bounded derived category of
 coherent sheaves on $X$ and the bounded derived category of finitely
 generated right modules over $A:=
 \End(\bigoplus_{0\leq i\leq \ell} L_i)$.
 \end{theorem}
  \begin{corollary}
   The algebra $\kk Q/R$ is a minimal noncommutative resolution of
   $\mathbb{A}^2_\kk/G$.
 \end{corollary}
 
 Wemyss~\cite{Wemyss} gives a complementary, algebraic approach to the
 two main results in this paper that is closer in spirit to
 noncommutative geometry than the approach adopted here.  The benefit
 of his approach is that the algebra $\kk Q/R$ is constructed without
 prior knowledge of the minimal resolution; the drawback is that the
 ideal of relations $R$ must be computed explicitly. We believe that both approaches have their
 merits.

\medskip
\noindent\textbf{Conventions}
  Write $\kk$ for an algebraically closed field of characteristic
   zero, $\kk^\times$ for the one-dimensional algebraic torus over
   $\kk$, and $\NN$ for the semigroup of nonnegative integers. We do not assume that toric
   varieties are normal. 
 
\medskip
\noindent\textbf{Acknowledgements.}  Special thanks to Alastair King and Michael Wemyss for pointing out errors in an earlier version. Thanks to Akira Ishii, Diane Maclagan and Greg Smith for useful discussions, and to Bill Crawley-Boevey for a helpful remark. 

\section{Background}
 \subsection{Minimal resolution by toric geometry}
 \label{sec:toric}
For a finite abelian subgroup $G\subset \GL(2,\kk)$ of order $r$, let $G^\vee=\Hom(G,\kk^\times)$ denote the character group and $\Irr(G)$ the set of equivalence classes of irreducible
 representations.  After killing quasireflections and changing
 coordinates if necessary, we may assume that $G$ is the cyclic group
 of order $r$ generated by the diagonal matrix $g =\textrm{diag}(\omega,
 \omega^a)$, where $\omega$ is a primitive $r$th root of unity and
 $\gcd(a,r)=1$. This action is said to be of type $\frac{1}{r}(1,a)$. The
 given representation of $G$ decomposes as $\rho_1\oplus \rho_2$,
 where $\rho_1(g)=\omega$ and $\rho_2(g)=\omega^a$.  The induced
 $G$-action on $\kk[x,y]$ satisfies $g\cdot x = \rho_1(g^{-1})x$
 and $g\cdot y = \rho_a(g^{-1})y$, and we obtain a $G^\vee$-grading of
 $\kk[x,y]$ via $\deg(x)=\rho_1$ and $\deg(y)=\rho_2$.

 To construct $\mathbb{A}_{\kk}^2/G$ and its
 minimal resolution by toric geometry, define
 $N:=\ZZ^2+\ZZ\cdot\frac{1}{r}(1,a)$ and $M:=\Hom_\ZZ(N,\ZZ)$. The
 nonnegative quadrant $\sigma$ in $N\otimes_{\ZZ}\QQ$ gives $\mathbb{A}_{\kk}^2/G =\Spec \kk[\sigma^\vee\cap M]$ with
 dense algebraic torus $T_M:= \Spec\kk[M]$.  For the minimal
 resolution, expand
 \begin{equation}
 \label{eqn:ctdfrac}
 \frac{r}{a} = c_{1} - \frac{1}{c_{2} - \frac{1}{\dots - \frac{1}{c_{\ell}}}},
 \end{equation}
 as a Jung-Hirzebruch continued fraction, giving integers
 $c_{1},\dots ,c_{\ell}\in \ZZ_{\geq 2}$. Define $(\beta_{\ell+1},\alpha_{\ell+1})=(0,r)$, $(\beta_\ell,\alpha_\ell)=(1,a)$, 
 and $(\beta_{i-1}, \alpha_{i-1}):= c_{i} (\beta_{i}, \alpha_{i}) - (\beta_{i+1}, \alpha_{i+1})$
 for $1\leq i \leq \ell$, which implies $(\beta_{0},\alpha_{0})=(r,0)$.  Define rays $\tau_0,\tau_1,\dots,\tau_{\ell+1}$ in $N\otimes_{\ZZ} \QQ$, where each $\tau_i$ has primitive generator $\frac{1}{r}(\beta_i,\alpha_i)\in N$.  The fan $\Sigma$ in $N\otimes_{\ZZ}\QQ$ obtained by subdividing the cone $\sigma$ by $\tau_1,\dots,\tau_\ell$ is the minimal resolution $f\colon X \to \mathbb{A}_{\kk}^2/G$. For $1\leq i\leq \ell$, let $D_{i}=X_{Star(\tau_i)}\cong \mathbb{P}^1$ denote the exceptional $T_M$-invariant divisor on $X$ with toric coordinates $[x^{\alpha_i} : y^{\beta_i}]$.  Write $\sigma_i$ for the cone generated by $\tau_i$ and $\tau_{i+1}$, so $X$ is covered by
 charts $U_i=\Spec\kk[M\cap\sigma_i^\vee]\cong
 \mathbb{A}^2_{\kk}$ for $0\leq i\leq \ell$, where
  \[
 U_i= \Spec
 \kk\left[\frac{x^{\alpha_{i+1}}}{y^{\beta_{i+1}}},\frac{y^{\beta_{i}}}{x^{\alpha_{i}}}\right].
 \]

  For a global description, write
 $\Sigma(1)=\{\tau_0,\tau_1,\dots,\tau_{\ell+1}\}$ for the set of rays
 in $\Sigma$,  and write $\NN^{\Sigma(1)}$ and $\ZZ^{\Sigma(1)}$ respectively for the semigroup and for the lattice generated by the $T_M$-invariant prime divisors.  Since $X$
 is smooth there is a short exact sequence
\begin{equation} 
  \label{eqn:degree}
  \begin{CD}   
    0 @>>> M @>>> \ZZ^{\Sigma(1)} @>\deg >> \Pic(X) @>>> 0, 
  \end{CD}
\end{equation}
where $\deg(D)=\mathscr{O}_X(D)$.  The total coordinate ring of $X$ is
the polynomial ring $\kk[x_0,\dots, x_{\ell +1}]$ obtained as the
semigroup algebra of $\NN^{\Sigma(1)}$. The degree map endows
$\kk[x_0,\dots,x_{\ell +1}]$ with a $\Pic(X)$-grading, and the
algebraic torus $\Hom(\Pic(X),\kk^\times)$ acts on
$\mathbb{A}^{\Sigma(1)}_\kk = \Spec\kk[x_0,\dots, x_{\ell+1}]$.  For
any line bundle $L\in \Pic(X)$, the $L$-graded piece
$\kk[x_0,\dots,x_{\ell+1}]_{L}$ is isomorphic to $H^0(X,L)$. In particular, if $L$ is
$f$-ample then $X$ is isomorphic to the GIT quotient
 \begin{equation}
 \label{eqn:toricGIT}
  \textstyle{\mathbb{A}_\kk^{\Sigma(1)}\git_{\! L} \Hom(\Pic(X),\kk^\times)}:=\Proj \Big{(}\bigoplus_{j \geq 0}\kk[x_0,\dots,x_{\ell+1}]_{L^{j}}\Big{)}.
  \end{equation}
This global description coincides with the quotient construction of $X$ from Cox~\cite{Cox}.

\subsection{Representations of bound quivers}
Let $Q$ be a finite connected quiver with vertex set $Q_0$, arrow set $Q_1$,
 and maps $\head, \tail \colon Q_1 \to Q_0$ indicating the vertices at
 the head and tail of each arrow.  The
 characteristic functions $\chi_{i} \colon Q_0 \to \ZZ$ for $i \in
 Q_0$ and $\chi_{a} \colon Q_1 \to \ZZ$ for $a \in Q_1$ form the
 standard integral bases of the vertex space $\ZZ^{Q_0}$ and the
 arrows space $\ZZ^{Q_1}$ respectively.  The incidence map $\inc
 \colon \ZZ^{Q_1} \to \ZZ^{Q_0}$ defined by
 $\inc(\chi_{a})=\chi_{\head(a)} - \chi_{\tail(a)}$ has image equal to
 the sublattice $\Wt(Q) \subset \ZZ^{Q_0}$ of functions $\theta \colon
 Q_0 \to \ZZ$ satisfying $\sum_{i \in Q_0} \theta_i = 0$. A nontrivial path in $Q$ is a
 sequence of arrows $p = a_1 \dotsb a_k$ with $\head(a_{j}) =
 \tail(a_{j+1})$ for $1 \leq j < k$.  We set $\tail(p) = \tail(a_{1}),
 \head(p) = \head(a_k)$ and $\supp(p)=\{a_1,\dots, a_k\}$. Each $i
 \in Q_0$ gives a trivial path $e_i$ where $\tail(e_i) = \head(e_i) =
 i$.  The path algebra $\kk Q$ is the $\kk$-algebra whose underlying
 $\kk$-vector space has a basis consisting of paths in $Q$, where the
 product of basis elements equals the basis element defined by
 concatenation of the paths if possible, or zero otherwise.

  A representation of a quiver $Q$ consists of a $\kk$-vector space $W_i$ for $i \in Q_0$ and a
 $\kk$-linear map $w_a \colon W_{\tail(a)} \to W_{\head(a)}$ for $a
 \in Q_1$. We often write $W$ as shorthand for $\bigl((W_i)_{i\in Q_0}, (w_a)_{a\in Q_1}\bigr)$. We consider only representations with
 $\dim_\kk(W_i)\leq 1$ for all $i\in Q_0$.  A map of representations
 $W\to W^\prime$ is a family $\psi_{i} \colon W_i^{\,} \to W_i^\prime$
 of $\kk$-linear maps for $i \in Q_0$ that are compatible with the
 structure maps, that is, $w_a^\prime \psi_{\tail(a)} =
 \psi_{\head(a)} w_a$ for all $a \in Q_1$.  With composition defined
 componentwise, we obtain the abelian category of representations of
 $Q$.   For any rational weight $\theta \in \Wt(Q)\otimes_{\ZZ}\QQ$, a
 representation $W$ satisfying $\dim_\kk(W_i)=1$ for all $i\in Q_0$ is
 $\theta$-stable if every proper, nonzero
 subrepresentation $W^\prime \subset W$ satisfies $\theta(W^\prime) :=
 \sum_{i \in \supp(W^\prime)} \theta_i > 0$, where $\supp(W') := \{ i
 \in Q_0 : W_i' \neq 0 \}$.  For $\theta$-semistability, replace $>$ with $\geq$.  

Let $R$ be a two-sided ideal in $\kk Q$ generated by differences of the form $p-q\in \kk Q$ where $p,q$ are paths with the same head and tail, each of which comprises at least two arrows. We do not assume that $R$ is admissible, so $\kk Q/R$ need not be of finite dimension over $\kk$. The pair $(Q,R)$ is an example of a \emph{bound quiver}, also known as a quiver with relations. For any representation $W=\bigl((W_i)_{i\in Q_0}, (w_a)_{a\in Q_1}\bigr)$ of $Q$ and for any nontrivial path $p = a_1 \dotsb
 a_k$ in $Q$, the evaluation of $W$ on $p$ is the  $\kk$-linear map
 $w_p \colon W_{\tail(p)} \to W_{\head(p)}$ defined by the composition
 $w_p = w_{a_1} \dotsb w_{a_k}$. A representation of the bound quiver $(Q,R)$ is a representation $W$ of $Q$ such that $w_p = w_{q}$ for all relations $p-q \in R$. The abelian category of finite-dimensional representations of $(Q,R)$ is equivalent to the category of finitely-generated $\kk Q/R$-modules.

\subsection{Moduli spaces of quiver representations}
Let $Q$ be a finite connected quiver. The space of representations $W$ of $Q$ for which $\dim_\kk(W_i)=1$ for all $i\in
 Q_0$ is
\[
 \textstyle{\mathbb{A}_{\kk}^{Q_1}:= \Spec \bigl( \kk[y_a : a\in Q_1] \bigr) \cong
\bigoplus\limits_{a \in Q_1} \Hom_{\kk}(W_{\tail(a)}, W_{\head(a)})}.
\]
The incidence map gives a $\Wt(Q)$-grading of the polynomial ring
$\kk[y_a : a\in Q_1]$ that induces an action of the algebraic torus
$T_Q:=\Hom(\Wt(Q),\kk^\times)$ on $\mathbb{A}^{Q_1}_\kk$, where $(t
\cdot w)_{a} = t_{\head(a)}^{\,} w_{a} t_{\tail(a)}^{-1}$.  For $\theta \in \Wt(Q)$, let $\kk[y_a : a\in Q_1]_{\theta}$
denote the $\theta$-graded piece. If every $\theta$-semistable representation of $Q$ is
 $\theta$-stable then King~\cite[Proposition~5.3]{King} proved that the GIT quotient
  \begin{equation}
 \label{eqn:quiverGIT}
\textstyle{\mathcal{M}_{\theta}(Q) := \mathbb{A}^{Q_1}_\kk\git_{\theta}T_Q = \Proj \Big{(}
\bigoplus_{j \geq 0} \kk[y_a : a\in Q_1]_{j
  \theta} \Big{)}}
 \end{equation}
 is the fine moduli space of isomorphism classes of $\theta$-stable
 representations of $Q$.  The $T_Q$-equivariant vector bundle
 $\bigoplus_{i \in Q_0} \mathscr{O}_{\mathbb{A}^{Q_1}}$ on $\mathbb{A}^{Q_1}_\kk$ descends to a
 tautological vector bundle $\bigoplus_{i \in Q_0} \mathscr{W}_i$ on $\mathcal{M}_\theta(Q)$. Our quivers will always have a distinguished vertex (denoted
 $0\in Q_0$ or $\rho_0\in Q_0^\prime$), and we normalise the universal
 family on $\mathcal{M}_\theta(Q)$ by identifying $T_Q$ with $\{
 (t_i)_{i\in Q_0} \in ( \kk^\times )^{Q_0} : t_0 = 1 \}$; this gives
 $\mathscr{W}_0\cong\mathscr{O}_{\mathcal{M}_\theta(Q)}$.

To construct moduli spaces of bound quiver representations, let $(Q,R)$ be a bound quiver where $R$ is generated by path differences $p-q\in \kk Q$. The map sending a path $p$ to the monomial $y_p:=\prod_{a\in
  \supp(p)} y_a\in \kk[y_a : a\in Q_1]$ enables us to define the ideal
 \begin{equation}
 \label{eqn:IR}
 I_R:= \big(y_p-y_q\in \kk[y_a : a\in Q_1] : p-q\in R\big)
 \end{equation}
 of relations in $\kk[y_a : a\in Q_1]$.  A point in $\mathbb{A}_\kk^{Q_1}$ corresponds to a
representation of $(Q,R)$ if and only it lies in the subscheme
$\mathbb{V}(I_R)$ cut out by $I_R$. The
binomial ideal $I_R$ is homogeneous with respect to the $\Wt(Q)$-grading, so
$\mathbb{V}(I_R)$ is invariant under the action of $T_Q$. If every $\vartheta$-semistable representation of $Q$ is $\vartheta$-stable then
\begin{equation}
\label{eqn:boundmoduli}
\mathcal{M}_\vartheta(Q,R):= \mathbb{V}(I_R)\git_\vartheta T_Q =
\textstyle{ \Proj \Big{(}\bigoplus_{j \geq 0} \big(\kk[y_a : a\in Q_1]/I_R)_{j \vartheta} \Big{)}}
\end{equation}
is the fine moduli space of isomorphism classes of $\vartheta$-stable
 representations of $(Q,R)$. The tautological bundle on $\mathcal{M}_\vartheta(Q,R)$ is obtained from that on $\mathcal{M}_\vartheta(Q)$ by restriction.
 
\section{Minimal resolution via multigraded linear series}

 \subsection{Quivers of sections and multigraded linear series}
 Let $\underline{\mathscr{E}}:=\big(\mathscr{E}_0,\mathscr{E}_1,\dots,\mathscr{E}_m\big)$
  be a sequence of distinct effective line bundles on the minimal resolution $X$ of $\mathbb{A}^2_\kk/G$, where $\mathscr{E}_0:=\mathscr{O}_X$ and $m>0$.  A $T_M$-invariant
 section $s \in H^0(X,\mathscr{E}_j^{\;} \otimes \mathscr{E}_{i}^{-1})=\Hom(\mathscr{E}_i,\mathscr{E}_j)$ is
 \emph{irreducible} if it does not factor through some $\mathscr{E}_k$ with $k\neq i,j$. The \emph{quiver of sections} of $\underline{\mathscr{E}}$ is the quiver $Q$ in which the vertex set $Q_0 = \{ 0,
 \dotsc,m \}$ corresponds to the line bundles in $\underline{\mathscr{E}}$, and
 where the arrows from $i$ to $j$ correspond to the irreducible sections in $H^0(X,\mathscr{E}_j^{\;} \otimes \mathscr{E}_{i}^{-1})$. For each arrow $a \in Q_1$, we write $\div(a) := \div(s) \in \NN^{\Sigma(1)}$ for the divisor of zeroes of the defining section $s \in H^0(X,\mathscr{E}_j^{\;} \otimes \mathscr{E}_{i}^{-1})$ and, more generally, for any path $p$ in $Q$ we call $\div(p) := \sum_{a\in \supp(p)} \div(a)$ the \emph{labelling divisor}.  The  \emph{ideal of relations} in the path algebra $\kk Q$ is the two-sided ideal
 \begin{equation}
 \label{eqn:relations}
 R = \big\langle p-q \in \kk Q : \head(p)=\head(q),
 \tail(p)=\tail(q), \div(p) =
 \div(q)\big\rangle
 \end{equation}
 Following Craw--Smith~\cite{CrawSmith} we call $(Q,R)$ the \emph{bound quiver of sections of} $\underline{\mathscr{E}}$. 
 
 \begin{lemma}
 \label{lem:algebra}
 Let $(Q,R)$ be the bound quiver of sections of
 $\underline{\mathscr{E}}$. The quiver $Q$ is connected and the quotient algebra $\kk Q/R$ is isomorphic to
 $\End\bigl( \bigoplus_{i\in Q_0}\mathscr{E}_i \bigr)$.
 \end{lemma}
 \begin{proof}
 The quiver is connected since $H^0(X,\mathscr{E}_i) \neq 0$ for $0 \leq i \leq m$. The algebra isomorphism follows as in the proof of \cite[Proposition~3.3]{CrawSmith}.
 \end{proof} 
 
\begin{remark}
 Since $X$ is projective over $\mathbb{A}^2_\kk/G$ we have $H^0(\mathscr{O}_X)\cong \kk[x,y]^G$, so directed cycles in $Q$ based at $i\in Q_0$ correspond to $G$-invariant monomials in $\kk[x,y]$. This makes the semiprojective situation different from the projective case; notably, the order $\mathscr{E}_1,\dots,\mathscr{E}_m$ is unimportant.
\end{remark} 
 
 The \emph{multigraded linear series} of $\underline{\mathscr{E}}$ is the variety $\vert\underline{\mathscr{E}}\vert := \mathbb{A}^{Q_1}_\kk\git_{\vartheta}T_Q$, where $Q$ is the quiver of sections of $\underline{\mathscr{E}}$ and where $\vartheta :=(-m,1,\dots,1)\in \Wt(Q)$ is a specially chosen weight. The next result extends the construction of Craw--Smith~\cite[Proposition~3.9]{CrawSmith} to the semiprojective setting.   
  
 \begin{proposition}
\label{prop:generic}
The multigraded linear series $\vert\underline{\mathscr{E}}\vert$ is isomorphic to the fine moduli space $\mathcal{M}_\vartheta(Q)$. This smooth toric variety is projective over the affine quotient $\mathbb{A}_\kk^{Q_1}/ T_Q$, and it is obtained as the geometric quotient of $\mathbb{A}_\kk^{Q_1}\setminus \mathbb{V}(B_Q)$ by $T_Q$ where 
 \[
 B_Q = \Bigg( \prod_{a\in \supp(T)} y_a : T \text{ is a spanning tree in } Q \text{ with root at }0\in Q_0\Bigg).
  \]
\end{proposition}
 \begin{proof}  
To prove the first statement it is enough to prove that every $\vartheta$-semistable representation is $\vartheta$-stable.  Let $W$ be a $\vartheta$-semistable representation of $Q$. If $W^\prime\subset W$ is a proper nonzero subrepresentation then $\sum_{i\in \supp(W^\prime)} \vartheta_i \geq 0$. Since $\vartheta_i=1$ for $i\neq 0$ and $\vartheta_0=-m$ we have $W_0^\prime = 0$ and $\vartheta(W^\prime)>0$, so $W$ is $\vartheta$-stable as required. A representation $W=\bigl((W_i)_{i\in Q_0}, (w_a)_{a\in Q_1}\bigr)$ is $\vartheta$-stable if and only if every subrepresentation $W^\prime\subset W$ with $W_0^\prime\neq 0$ has $W_i^\prime\neq 0$ for all $i\neq 0$, which holds if and only if there is a spanning tree $T$ with root at 0 such that $\prod_{a\in \supp(T)} w_a\neq 0$. The ideal $B_Q$ therefore cuts out the $\vartheta$-unstable locus. That the toric variety $\mathbb{A}^{Q_1}_\kk\git_{\!  \vartheta} T_Q$ is smooth follows from the fact that the incidence map of $Q$ is totally unimodular as in \cite[Proposition~3.8]{CrawSmith}, and variation of GIT quotient $\vartheta\mapsto 0$ gives a projective morphism from $\mathbb{A}^{Q_1}_\kk\git_{\!  \vartheta} T_Q$ to $\mathbb{A}_\kk^{Q_1}/T_Q$.
  \end{proof}

Since $\mathcal{M}_\vartheta(Q)$ is a fine moduli space, it follows that the multigraded linear series $\vert\underline{\mathscr{E}}\vert$ carries tautological line bundles $\mathscr{W}_0, \dots,\mathscr{W}_m$ where $\mathscr{W}_0\cong \mathscr{O}_{\vert\underline{\mathscr{E}}\vert}$.
 
\subsection{A morphism to the multilinear series}
We choose once and for all a preferred sequence of line bundles on $X$, namely, the sequence 
\begin{equation}
\label{eqn:sequence}
\underline{\mathscr{L}}=\big(\mathscr{O}_X, L_1,\dots, L_\ell\big),
 \end{equation}
  such that $\deg(L_i\vert_{D_j}) = \delta_{ij}$ for $1\leq i,j\leq \ell$. These line bundles are nef and hence globally generated; in explicit toric coordinates,  for $0\leq j\leq \ell$, the free $\mathscr{O}_{U_j}$-module $L_i\vert_{U_j}$ of rank one is generated by $x^{\alpha_i}$ if $i\leq j$, and $y^{\beta_i}$ if $i > j$.   
  
 The multigraded linear series of $\underline{\mathscr{L}}$ is the fine moduli space $\vert\underline{\mathscr{L}}\vert:=\mathcal{M}_\vartheta(Q)$, where $Q$ is the quiver of sections of $\underline{\mathscr{L}}$ and $\vartheta = (-\ell,1,\dots,1)\in \Wt(Q)$.   The
 map $\Phi_Q\colon \kk[y_a : a\in Q_1]\to \kk[\NN^{\Sigma(1)}]$
 defined by $\Phi_Q(y_a) = x^{\div(a)}$ induces a morphism
 $(\Phi_Q)^*\colon \mathbb{A}^{\Sigma(1)}_\kk \to
 \mathbb{A}^{Q_1}_\kk$ that is equivariant with respect to the actions
 of $\Hom(\Pic(X),\kk^\times)$ and
 $T_Q=\Hom(\Wt(Q),\kk^\times)$. Since each $L_i$ is globally generated, \cite[Section~4]{CrawSmith} shows that
 $(\Phi_Q)^*$ descends to give a morphism $\varphi_{\vert
 \underline{\mathscr{L}}\vert}\colon X\longrightarrow  \vert\underline{\mathscr{L}}\vert$.
 To describe the image, write $\pi := (\inc,\div)
 \colon \ZZ^{Q_{1}} \to \Wt(Q) \oplus \NN^{\Sigma(1)}$ for the
 $\ZZ$-linear map sending $\chi_{a}$ to $\bigl( \chi_{\head(a)} -
 \chi_{\tail(a)}, \div(a) \bigr)$ for $a\in Q_1$. Write $\NN(Q)$ for
 the image under $\pi$ of the subsemigroup $\NN^{Q_1}$ generated by
 $\chi_a$ for $a\in Q_1$, and write $\kk[\NN(Q)]$ for its semigroup
 algebra. The projections $\pi_1\colon \NN(Q)\to
 \Wt(Q)$ and $\pi_2\colon \NN(Q)\to \NN^{\Sigma(1)}$ fit in to the
 commutative diagram
\begin{equation}
 \label{eqn:diagram}
\xymatrix{\NN(Q)\ar[r]^*{\text{\footnotesize{$\pi_1$}}}\ar[d]^*{\text{\footnotesize{$\pi_2$}}} & \Wt(Q)\ar[d]^*{\text{\footnotesize{$\pic$}}} \\
 \NN^{\Sigma(1)} \ar[r]^*{\text{\footnotesize{$\deg$}}} & \Pic(X)}
\end{equation}
where $\pic(\sum_{i\in Q_0} \theta_i \chi_{i}):=\bigotimes_{i\in Q_0}
L_{i}^{\theta_{i}}$. Since $\kk[y_a : a\in Q_1]$ is the semigroup
algebra of $\NN^{Q_1}$, the map $\pi$ induces a surjective maps of
$\kk$-algebras $\pi_*\colon \kk[y_a: a\in Q_1]\to \kk[\NN(Q)]$ with
kernel 
 \begin{equation}
 \label{eqn:IQ}
 I_Q := \big(y^u-y^v \in \kk[y_a : a\in Q_1] ; u-v\in \Ker(\pi)\big).
\end{equation}
 The incidence map factors through
$\NN(Q)$ to give the map $\pi_1$, so the action of $T_Q$ on $\mathbb{A}_\kk^{Q_1}$ restricts
to an action on the toric variety $\mathbb{V}(I_Q) = \Spec
\kk[\NN(Q)]$ cut out by the toric ideal $I_Q$.

 \begin{theorem}
 \label{thm:closedimmersion}
 Let $Q$ denote the quiver of sections of the sequence of bundles $\underline{\mathscr{L}}$ from \eqref{eqn:sequence}. The morphism $\varphi_{\vert \underline{\mathscr{L}}\vert}\colon X\to
 \vert\underline{\mathscr{L}}\vert$ is a closed immersion with image $\mathbb{V}(I_Q) \git_{\! \vartheta}T_Q$, and the tautological line bundles $\mathscr{W}_0, \dots,\mathscr{W}_\ell$ on $\vert\underline{\mathscr{L}}\vert$ satisfy $\varphi_{\vert \underline{\mathscr{L}}\vert}^*(\mathscr{W}_i) = L_i$ for $0\leq i\leq \ell$.
 \end{theorem}
 \begin{proof}
   The toric ideal $I_Q$
   is the $\Wt(Q)$-homogeneous part of $\Ker(\Phi_Q)$.  Thus, just as
   $\Ker(\Phi_Q)$ cuts out the image of $(\Phi_Q)^*$, so $I_Q$ cuts
   out the image of $\varphi_{\vert \mathcal{L}\vert}$ after passing
   to the quotient by the action of $T_Q$. Following
   \cite[Theorem~1.1]{CrawSmith}, the image of $\varphi_{\vert
     \mathcal{L}\vert}$ is the GIT quotient $\mathbb{V}(I_Q) \git_{\!
     \vartheta}T_Q$. 
     
     To see that $\varphi_{\vert \underline{\mathscr{L}}\vert}$ is
   a closed immersion, set $\mathcal{V}:= \NN(Q)\cap
   \pi_1^{-1}(\vartheta)$ and $\mathcal{W}:=\NN^{\Sigma(1)}\cap
   \deg^{-1}(L)$ for the very ample line bundle $L=\bigotimes_{i=1}^\ell
   L_i$.  Since $L=\pic(\vartheta)$, we have
   $\pi_2(\mathcal{V})\subseteq \mathcal{W}$ and hence $\Phi_Q(\kk[y_a
   : a\in Q_1]_\vartheta)\subseteq \kk[x_0,\dots,x_{\ell+1}]_L \cong
   H^0(X,L)$.  According to the proof of
   \cite[Theorem~4.9]{CrawSmith}, the morphism $\varphi_{\vert\underline{\mathscr{L}}\vert}$ is a closed immersion if and only if the
   linear series $\Phi_Q(\kk[y_a : a\in Q_1]_\vartheta )$ defines a
   closed immersion.  As a result, we must show that the $\kk$-vector space
   basis $\pi_2(\mathcal{V})$ of $\Phi_Q(\kk[y_a : a\in
   Q_1]_\vartheta)$ contains the vertices $\{u(\sigma_j) : 0\leq j\leq
   \ell\}$ of the polyhedron $P_L$ obtained as the convex hull of the
   set $\mathcal{W}$ and, in addition, that the semigroup $M\cap
   \sigma_j^\vee$ is generated by $\{u-u(\sigma_j) : u\in
   \pi_2(\mathcal{V})\}$ for each $0\leq j\leq \ell$.    First, identify $\NN(Q)\cap \pi_1^{-1}(\chi_i-\chi_0)$ with the monomial basis of $H^0(X,L_i)$ for each $0\leq i\leq \ell$, and write
   $v_j(i)\in \NN(Q)\cap \pi_1^{-1}(\chi_i-\chi_0)$ for the element
   $x^{\alpha_i}\in H^0(X,L_i)$ if $0\leq i\leq j$, and
   $y^{\beta_i}\in H^0(X,L_i)$ if $j<i\leq \ell$. Since  $L\vert_{U_j}$ is generated as an
   $\mathscr{O}_{U_j}$-module by
   $s_j:=\prod_{0\leq i\leq j} x^{\alpha_i}\cdot
   \prod_{j<i\leq\ell}y^{\beta_i}$, the element
   $v_j:=\sum_{i=0}^{\ell}v_j(i)\in \mathcal{V}$ satisfies
   $\pi_2(v_j)=\div(s_j)=u(\sigma_j)$, so the vertices of $P_L$ lie in
   $\pi_2(\mathcal{V})$.  For the second part, write $v_{\alpha_i},
   v_{\beta_i}\in\NN(Q)\cap \pi_1^{-1}(\chi_i-\chi_0)$ for the
   elements $x^{\alpha_i},y^{\beta_i}\in
   H^0(X,L_i)$ respectively, where $i=j,j+1$.  Define elements of $\mathcal{V}$ by
   setting $v_j^+:= v_j+(v_{\beta_j}-v_{\alpha_j})$ and $v_j^-:=
   v_j+(v_{\alpha_{j+1}}-v_{\beta_{j+1}})$. We obtain $\pi_2(v_j^+) =
   \div(s_j\cdot(y^{\beta_j}/x^{\alpha_j}))$ and $\pi_2(v_j^-) =
   \div(s_j\cdot(x^{\alpha_{j+1}}/y^{\beta_{j+1}}))$, so for each $j$
   the differences $v_j^+-v_j,v_j^{-}-v_j\in \mathcal{V}-v_j$ map
   under $\pi_2$ to the generators $\div(y^{\beta_j}/x^{\alpha_j}),
   \div(x^{\alpha_{j+1}}/y^{\beta_{j+1}})$ of $M\cap \sigma_j^\vee$.
   Thus, the morphism $\varphi_{\vert \underline{\mathscr{L}}\vert}$ is indeed a closed immersion. For the statement about the tautological bundles, note that the proof of \cite[Theorem~4.15]{CrawSmith} applies verbatim in this case.
 \end{proof}

 \begin{corollary}
 The minimal resolution $f\colon
 X\to \mathbb{A}^2_\kk/G$ coincides with the projective morphism
 $\mathbb{V}(I_Q)\git_\vartheta T_Q\to \mathbb{V}(I_Q)/ T_Q$
 obtained by variation of GIT quotient.
 \end{corollary}
 \begin{proof}
 This is similar to the proof of \cite[Proposition~4.1, Theorem~4.3\one]{CMT1}.
 \end{proof}
 
 \begin{example}
 \label{ex:712}
  For the action of type $\frac{1}{7}(1,2)$, the rational curves
 $D_1$ and $D_2$ in $X$ have toric coordinates $[x:y^4]$ and $[x^2:y]$
 respectively, as is evident from the fan of $X$ shown in
 Figure~\ref{fig:712}(a).  The sequence is $\underline{\mathscr{L}} =
 (\mathscr{O}_X, L_1, L_2)$, where: $L_1\vert_{U_j}$ is generated by
 $x$ for $j=1,2$ and $y^4$ for $j=0$; while $L_2\vert_{U_j}$ is
 generated by $x^2$ for $j=2$ and $y$ for $j=0,1$. 
 Since $L_1$ has degree 1 on $D_1$ and degree 0 on
 $D_2$, we have $L_1=\mathscr{O}_X(D_0)$ and, similarly, $L_2 = \mathscr{O}_X(D_3)$. 
   \begin{figure}[!ht]
    \centering
      \subfigure[Toric fan]{
       \psset{unit=1cm}
     \begin{pspicture}(0,-0.5)(2.5,2.75)
       \psline{*-*}(0,0)(2.5,0)
       \psline{*-}(0,0)(1.25,2.5)
       \psline{*-}(0,0)(2.5,0.625)
       \psline{*-*}(0,0)(0,2.5)
       \psdot(0.35714,0.71429)
       \psdot(0.71429,1.42857)
       \psdot(1.07151,2.14287)
       \psdot(1.42857,0.357145)
       \psdot(1.78571,1.071435)
       \psdot(2.14285,1.785725)
       \rput(0.5,2){$\sigma_2$}
       \rput(1.2,0.9){$\sigma_1$}
       \rput(2.1,0.2){$\sigma_0$}
       \rput(1.4,2.3){$\tau_2$}
       \rput(2.4,0.8){$\tau_1$}
       \end{pspicture}
       }
      \qquad  \qquad
      \subfigure[Quiver of sections]{
        \psset{unit=1.3cm}
        \begin{pspicture}(-0.5,-0.5)(2,2.75)
        \cnodeput(0,0){A}{1}
        \cnodeput(2.5,0){B}{2} 
        \cnodeput(0,2.5){C}{0}
        \psset{nodesep=0pt}
        \nccurve[angleA=10,angleB=170]{->}{A}{B}\lput*{:U}{$x_0x_1$}
        \nccurve[angleA=345,angleB=195]{<-}{A}{B}\lput*{:U}{$x_2x_3^3$}
        \nccurve[angleA=285,angleB=75]{->}{C}{A}\lput*{:180}{$x_0$}
        \nccurve[angleA=255,angleB=105]{<-}{C}{A}\lput*{:180}{$x_1x_2x_3^3$}
        \nccurve[angleA=300,angleB=150]{<-}{C}{B}\lput*{:U}{$x_0^5x_1^3x_2$}
        \nccurve[angleA=320,angleB=130]{->}{C}{B}\lput*{:U}{$x_3$}
        \nccurve[angleA=345,angleB=105]{<-}{C}{B}\lput*{:U}{$x_0^3x_1^2x_2x_3$}
        \nccurve[angleA=15,angleB=75]{<-}{C}{B}\lput*{:U}{$x_0x_1x_2x_3^2$}
        \end{pspicture}}
      \qquad \qquad  
      \subfigure[Listing the arrows]{
        \psset{unit=1.3cm}
\begin{pspicture}(-0.5,-0.5)(3,2.5)
          \cnodeput(0,0){A}{1}
          \cnodeput(2.5,0){B}{2} 
          \cnodeput(0,2.5){C}{0}
          \psset{nodesep=0pt}
          \nccurve[angleA=10,angleB=170]{->}{A}{B}\lput*{:U}{$a_3$}
          \nccurve[angleA=345,angleB=195]{<-}{A}{B} \lput*{:U}{$a_4$}
          \nccurve[angleA=285,angleB=75]{->}{C}{A} \lput*{:90}{$a_1$}
          \nccurve[angleA=255,angleB=105]{<-}{C}{A} \lput*{:90}{$a_2$}
          \nccurve[angleA=300,angleB=150]{<-}{C}{B} \lput*{:U}{$a_5$}
          \nccurve[angleA=320,angleB=130]{->}{C}{B} \lput*{:U}{$a_6$}
          \nccurve[angleA=345,angleB=105]{<-}{C}{B} \lput*{:U}{$a_7$}
          \nccurve[angleA=15,angleB=75]{<-}{C}{B} \lput*{:U}{$a_8$}
        \end{pspicture}     
        }
    \caption{A quiver of sections for the action of type $\frac{1}{7}(1,2)$}
  \label{fig:712}
  \end{figure}
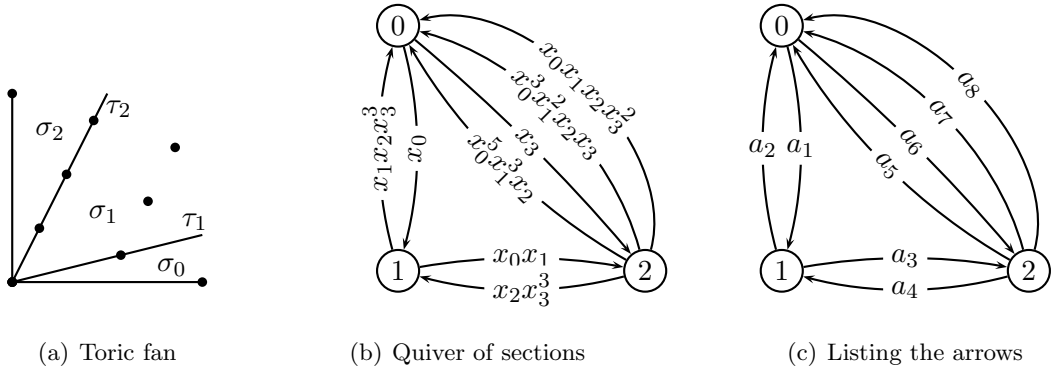
  The quiver of sections is shown in Figure~\ref{fig:712}(b), where each labelling divisor is recorded as a monomial in the Cox ring $\kk[x_0,x_1,x_2,x_3]$.
  One approach to drawing the quiver of sections of $\underline{\mathscr{L}}$ is to first
  draw the quiver of sections of $\mathscr{O}_X$; that is, one vertex
  and five loops, each labelled with a $\kk$-algebra generator of
  $H^0(\mathscr{O}_X)$. Adding a vertex for $L_1$
  and then $L_2$ causes the loops to decompose into
  paths according to the decomposition of the labelling divisor. In this case, the semigroup $\NN(Q)$ is generated by the
   columns of the matrix \renewcommand{\arraystretch}{0.8}
   \renewcommand{\arraycolsep}{2pt}
  \[
    \left[ \text{\footnotesize $\begin{array}{rrrrrrrr}
     -1 &  1 &  0 &  0 &  1 &  -1 &  1 &  1 \\
      1 & -1 & -1 &  1 &  0 &  0 &  0 &  0 \\
      0 &  0 &  1 & -1 &  -1 & 1 & -1 & -1 \\     
      1 &  0 &  1 &  0 &  5 &  0 &  3 &  1 \\
      0 &  1 &  1 &  0 &  3 &  0 &  2 &  1 \\
      0 &  1 &  0 &  1 &  1 &  0 &  1 &  1 \\
      0 &  3 &  0 &  3 &  0 &  1 &  1 &  2
      \end{array}$} \right] ,
  \]
  where the $i$-th column corresponds to $a_i$ for $1 \leq i \leq 8$.
  One computes using Macaulay2~\cite{M2} that
  \[
 I_Q= \left(
      \begin{array}{c}
      y_7^2-y_5y_8,\: y_3y_4-y_6y_8,\:  y_1y_2-y_6y_8,\: y_3y_7y_8-y_2y_5,\: y_1y_7y_8-y_4y_5, \\
y_1y_3y_7-y_5y_6,\: y_3y_8^2-y_2y_7,\: y_1y_8^2-y_4y_7,\: y_1y_3y_8-y_6y_7
  \end{array}
    \right)
   \]
 is the toric ideal that cuts out the image of $X$ under the morphism  $\varphi_{\vert \underline{\mathscr{L}}\vert}\colon X\to
 \vert\underline{\mathscr{L}}\vert$. 
 \end{example}
 
\begin{remark} 
\label{rem:IRinIQ}
It is easy to see that the ideal $I_R$ arising from $(Q,R)$ as in \eqref{eqn:IR} is contained in $I_Q$, so the image $\mathbb{V}(I_Q)\git_\vartheta T_Q$ of the closed immersion $\varphi_{\vert \underline{\mathscr{L}}\vert}\colon X\to \vert\underline{\mathscr{L}}\vert$ is a subvariety of the fine moduli space $\mathcal{M}_\vartheta(Q,R)$ of $\vartheta$-stable representations of the bound quiver $(Q,R)$.  
\end{remark}

\section{The bound Special McKay quiver}
\label{sec:specialMcKay}
\subsection{The bound McKay quiver}
 The \emph{McKay
   quiver} of the $G$-action of type $\frac{1}{r}(1,a)$ is the quiver $Q^\prime$ with
 vertex set $Q^\prime_0=\Irr(G)$ and arrow set $Q^\prime_1=\{a_1^\rho,
 a_2^\rho : \rho\in \Irr(G)\}$, where arrow $a_i^\rho$ goes from
 vertex $\rho\rho_i:=\rho\otimes\rho_i$ to vertex $\rho$ for all
 $\rho\in \Irr(G)$ and $i=1,2$.  The \emph{label} of each arrow $a_1^\rho$ is the monomial $\mon(a_1^\rho)=x$, and similarly, the label on $a_2^\rho$ is $\mon(a_2^\rho)=y$. More generally, the label of a path $p^\prime$ in $Q^\prime$ is the product $\mon(p^\prime) = \prod_{a^\prime\in \supp(p^\prime)} \mon(a^\prime)$. Let $\kk Q^\prime$ denote the path algebra of the quiver
 $Q^\prime$, and consider the two-sided ideal in $\kk Q^\prime$ given by
 \begin{equation}
 \label{eqn:McKayrelations}
 R^\prime:=\big\langle a_2^{\rho\rho_1}a_1^{\rho} -
 a_1^{\rho\rho_2}a_2^{\rho} : \rho\in\Irr(G)\big\rangle
 \end{equation} 
 Equivalently, $R^\prime$ is generated by path differences $p^\prime-q^\prime\in \kk Q$ for which $p^\prime, q^\prime$ have the same head, tail and labelling monomial. The pair $(Q^\prime,R^\prime)$ is the McKay
 quiver with relations, or equivalently, the \emph{bound McKay quiver}, of the subgroup
 $G\subset \GL(2,\kk)$. 
  
 The McKay correspondence provides a strong link between the bound McKay quiver and the minimal resolution $X$ of $\mathbb{A}^2_\kk/G$.  To state the result, let $\{\chi_\rho : \rho\in \Irr(G)\}$ denote the standard basis of
 the vertex space $\ZZ^{\Irr(G)}$ of $Q^\prime$, where $\rho_0\in\Irr(G)$ denotes the trivial
 representation.
 
 \begin{lemma}
 \label{lem:boundMcKay}
 Let  $(Q^\prime,R^\prime)$  denote the bound McKay quiver and set  $\vartheta^\prime:=\sum_{\rho\in
 \Irr(G)}(\chi_\rho-\chi_{\rho_0})$. The fine moduli space
 $\mathcal{M}_{\vartheta^\prime}(Q^\prime,R^\prime)$ of
 $\vartheta^\prime$-stable representations is isomorphic to $X$. Moreover, if $\bigoplus_{\rho\in \Irr(G)} \mathscr{W}_{\rho}$ denotes the tautological bundle on $\mathcal{M}_{\vartheta^\prime}(Q^\prime,R^\prime)$, then the following are isomorphic:
 \begin{enumerate}
  \item[\one] the skew group algebra $\kk[x,y]*G$;
  \item[\two] the quotient algebra $\kk Q^\prime/R^\prime$; and
   \item[\three] the endomorphism algebra $\End\big{(}\bigoplus_{\rho\in \Irr(G)}
 \mathscr{W}_{\rho}\big{)}$.
 \end{enumerate}
 \end{lemma}
 \begin{proof}
 The fine moduli space
 $\mathcal{M}_{\vartheta^\prime}(Q^\prime,R^\prime)$ is isomorphic to
 the $G$-Hilbert scheme $\ghilb(\mathbb{A}^2_\kk)$, so the first statement follows from Kidoh~\cite{Kidoh}. For the second, the isomorphism between \one\ and \two\ is well known (see \cite{CMT2}). For the isomorphism between \two\ and \three, let $\Gamma\subset \SL(3,\kk)$ be the subgroup of type
   $\frac{1}{r}(1,a,r-a-1)$ whose McKay quiver $Q_\Gamma$ can be
   constructed from that of $G\subset \GL(2,\kk)$ by adding an arrow
   labelled $z$ from $\rho$ to $\rho\rho_1\rho_2$ for all $\rho\in
   \Irr(\Gamma)=\Irr(G)$, and the corresponding ideal of relations
   $R_\Gamma\subset \kk Q_\Gamma$ includes additional generators to
   ensure that arrows labelled $z$ commute with those labelled $x$ and
   those labelled $y$. The McKay correspondence for $\Gamma$ by Bridgeland--King--Reid~\cite{BKR} 
   gives $\kk
   Q_\Gamma/R_\Gamma\cong
   \End_{\mathscr{O}_Y}\big{(}\bigoplus_{\rho\in \Irr(\Gamma)}
   \mathscr{R}_{\rho}\big{)}$, where $\bigoplus_{\rho\in \Irr(\Gamma)}
   \mathscr{R}_{\rho}$ is the tautological bundle on
   $Y=\Gammahilb(\mathbb{A}^3_\kk)$. Since the $G$-action on
   $\mathbb{A}^2_\kk$ is obtained from the $\Gamma$-action on
   $\mathbb{A}^3_\kk$ by setting $z=0$, functoriality of the
   $G$-Hilbert scheme implies that $X=\ghilb(\mathbb{A}^2_\kk)$ is a
   divisor of $Y=\Gammahilb(\mathbb{A}^3_\kk)$, and the
   tautological bundle $\bigoplus_{\rho\in \Irr(G)} \mathscr{W}_\rho$
   on $X$ is obtained from that on $Y$ by setting $z=0$ throughout.
   The result follows since $\kk Q/R$ is obtained from $\kk
   Q_\Gamma/R_\Gamma$ by setting to zero all arrows labelled $z$.
 \end{proof}
 
 The fine moduli space
 $\mathcal{M}_{\vartheta^\prime}(Q^\prime,R^\prime)$ represents the
 same functor as the $G$-Hilbert scheme, so the tautological bundles
 $\{\mathscr{W}_{\rho}: \rho\in \Irr(G)\}$ on
 $\mathcal{M}_{\vartheta^\prime}(Q^\prime,R^\prime)$ coincide with
 those on $\ghilb(\mathbb{A}^2_\kk)$. To characterise these bundles,
 Kidoh~\cite[Theorem~5.1]{Kidoh} observed that for $0\leq j\leq
 \ell$ and for the toric chart $U_j= \Spec
 \kk[x^{\alpha_{j+1}}/y^{\beta_{j+1}},y^{\beta_{j}}/x^{\alpha_{j}}]$ on
 $X$, the unique standard monomial of the ideal $I_j=(x^{\alpha_{j+1}},
 y^{\beta_{j}},x^{\alpha_{j+1}-\alpha_{j}}y^{\beta_{j}-\beta_{j+1}})$ in
 degree $\rho\in \Irr(G)$ generates the rank one
 $\mathscr{O}_{U_j}$-module $\mathscr{W}_{\rho^*}\vert_{U_j}$, where
 $\rho^*$ is the contragradient representation (recall from
 Section~\ref{sec:toric} that $\rho:=\deg(x^iy^j)$ satisfies
 $\rho(g)=\omega^{i+aj}$).  In particular, each $\mathscr{W}_\rho$ is
 nef.
 
  \begin{remark}
 \label{rem:contragradient}
  We emphasise that the monomials defining sections of the tautological bundle
   $\mathscr{W}_{\rho^*}$ on
   $\mathcal{M}_{\vartheta^\prime}(Q^\prime,R^\prime)$ indexed by
   vertex $\rho^*\in \Irr(G)$ have degree $\rho\in\Irr(G)$. 
 \end{remark}

 %The following result strengthens Lemma~\ref{lem:boundMcKay}\two.

 \begin{proposition}
 \label{prop:McKayquiverofsections}
   The bound McKay quiver $(Q^\prime, R^\prime)$ is the bound
   quiver of sections for the sequence of tautological line bundles
   $\underline{\mathscr{L}}^\prime:=(\mathscr{W}_{\rho} : \rho\in \Irr(G))$ on
   $\mathcal{M}_{\vartheta^\prime}(Q^\prime, R^\prime)\cong X$.
 \end{proposition}
 \begin{proof}
   Each tautological line bundle $\mathscr{W}_\rho$ is nef and
   $\mathscr{W}_{\rho_0}\cong\mathscr{O}_X$, so the bound
   quiver of sections $(Q_{\underline{\mathscr{L}}^\prime},
   R_{\underline{\mathscr{L}}^\prime})$ for $\underline{\mathscr{L}}^\prime$ is well defined.
   The vertex sets of $Q^\prime$ and $Q_{\underline{\mathscr{L}}^\prime}$
   coincide.  For any closed point $[W]\in
   \mathcal{M}_{\vartheta^\prime}(Q^\prime,R^\prime)$, the fibre
   $\mathscr{W}_\rho\vert_{[W]}$ is the $\kk$-vector space $W_\rho$
   encoded in the representation of $Q^\prime$ parametrised by $[W]$.
   As the point $[W]$ varies in
   $\mathcal{M}_{\vartheta^\prime}(Q^\prime,R^\prime)$, the maps
   $w_i^\rho \colon W_{\rho\rho_i}\to W_\rho$ arising from arrows
   $a_i^\rho\in Q^\prime_1$ determine sections $s_i^\rho \in
   \Hom(\mathscr{W}_{\rho\rho_i},\mathscr{W}_\rho)=H^0(\mathscr{W}_\rho
   \otimes \mathscr{W}_{\rho\rho_i}^{-1})$. These sections are
   irreducible because $R^\prime$ lies in the ideal of $\kk Q^\prime$
   generated by paths of length two.  Each arrow in $Q^\prime$
   therefore determines an arrow in $Q_{\mathcal{L}^\prime}$ with the
   same head and tail, so $Q^\prime$ is a subquiver of
   $Q_{\mathcal{L^\prime}}$.  For $a_2^{\rho\rho_1}a_1^{\rho} -
   a_1^{\rho\rho_2}a_2^{\rho}\in R^\prime$, the relation
   $w_2^{\rho\rho_1}w_1^\rho= w_1^{\rho\rho_2}w_2^{\rho}\colon
   W_{\rho\rho_1\rho_2}\to W_\rho$ determines a relation
   $s_2^{\rho\rho_1}s_1^\rho= s_1^{\rho\rho_2}s_2^{\rho}\colon
   \mathscr{W}_{\rho\rho_1\rho_2}\to \mathscr{W}_\rho$ between
   sections. Thus, $R^\prime$ is a subset of $R_{\mathcal{L^\prime}}$
   under the inclusion of $\kk Q^\prime$ as a subalgebra of $\kk
   Q_{\mathcal{L^\prime}}$. Thus far we have $Q^\prime\subseteq Q_{\mathcal{L^\prime}}$ and $R^\prime\subseteq R_{\mathcal{L^\prime}}$. Lemma~\ref{lem:algebra} and Lemma~\ref{lem:boundMcKay} together give an isomorphism $\kk
   Q^\prime/R^\prime\cong\kk
   Q_{\mathcal{L}^\prime}/R_{\mathcal{L}^\prime}$. If there exists an
   arrow $a\in Q_{\mathcal{L}^\prime}\setminus Q^\prime$ then this
   algebra isomorphism forces a relation in $R_{\mathcal{L}^\prime}$
   that contains a term of the form $\lambda a + f$ for some
   $\lambda\in \kk$ and $f\in \kk Q_{\mathcal{L}^\prime}$. However,
   $R_{\mathcal{L}^\prime}$ is an ideal of relations in a quiver of
   sections and hence every generating relation is a combination of
   paths of length at least two, so $Q^\prime= Q_{\mathcal{L}^\prime}$
   after all. In particular, $\kk Q^\prime = \kk
   Q_{\mathcal{L^\prime}}$. The isomorphism $\kk
   Q^\prime/R^\prime\cong\kk
   Q_{\mathcal{L}^\prime}/R_{\mathcal{L}^\prime}$ now forces the
   inclusion $R^\prime\subseteq R_{\mathcal{L}^\prime}$ to be equality
   as required.
 \end{proof}
 
 \begin{corollary}
 \label{coro:labels} 
 For paths $p^\prime, q^\prime$ in $Q^\prime$ with the same head and tail, the following are equivalent:
 \begin{enumerate}
 \item[\one] $\div(p^\prime)=\div(q^\prime)\in \kk[x_0,\dots,x_{\ell+1}]$.
 \item[\two] $p^\prime - q^\prime \in R^\prime$;
 \item[\three] $\mon(p^\prime)=\mon(q^\prime)\in \kk[x,y]$;
 \end{enumerate}
 \end{corollary}
    
   \begin{example}
 \label{ex:712two}
 Figure~\ref{fig:McKayquiver}(a) shows the McKay quiver for the action of type $\frac{1}{7}(1,2)$ with labels $x$ and $y$, while Figure~\ref{fig:McKayquiver}(b) gives the labels on $Q^\prime$ as a quiver of sections.   Replace $x_0$ by $x$, $x_3$ by $y$ and set $x_1=x_2=1$ to recover the labels in Figure~\ref{fig:McKayquiver}(a) from those in Figure~\ref{fig:McKayquiver}(b).
 \begin{figure}[!ht]
    \centering
    \mbox{
      \subfigure[]{
        \psset{unit=1.3cm}
        \begin{pspicture}(0.5,-1.2)(4.5,3.4)
        \cnodeput(2.5,3){A}{$\rho_0$}
        \cnodeput(1.2,2.3){B}{$\rho_1$} 
        \cnodeput(0.8,1){C}{$\rho_2$}
        \cnodeput(1.8,0){D}{$\rho_3$}
        \cnodeput(3.2,0){E}{$\rho_4$} 
        \cnodeput(4.2,1){F}{$\rho_5$}
        \cnodeput(3.8,2.3){G}{$\rho_6$}
        \psset{nodesep=0pt}
        \ncline{->}{B}{A}\lput*{:U}{$x$}
        \ncline{->}{C}{B}\lput*{:U}{$x$}
        \ncline{->}{D}{C}\lput*{:U}{$x$}
        \ncline{->}{E}{D}\lput*{:U}{$x$}
        \ncline{->}{F}{E}\lput*{:U}{$x$}
        \ncline{->}{G}{F}\lput*{:U}{$x$}
        \ncline{->}{A}{G}\lput*{:U}{$x$}
        \ncline{->}{C}{A}\lput*{:U}{$y$}
        \ncline{->}{D}{B}\lput*{:270}{$y$}
        \ncline{->}{E}{C}\lput*{:200}{$y$}
        \ncline{->}{F}{D}\lput*{:168}{$y$}
        \ncline{->}{G}{E}\lput*{:110}{$y$}
        \ncline{->}{A}{F}\lput*{:U}{$y$}
        \ncline{->}{B}{G}\lput*{:U}{$y$}
        \end{pspicture}}
      \qquad  \qquad
      \subfigure[]{
        \psset{unit=1.3cm} \begin{pspicture}(-0.5,-0.5)(5,4.4)
        \cnodeput(2.5,4.2){A}{$\mathscr{O}_X$}
        \cnodeput(0.6,3.3){B}{$\mathscr{W}_{\rho_6^*}$} 
        \cnodeput(0.1,1.6){C}{$\mathscr{W}_{\rho_5^*}$}
        \cnodeput(1.3,0){D}{$\mathscr{W}_{\rho_4^*}$}
        \cnodeput(3.7,0){E}{$\mathscr{W}_{\rho_3^*}$} 
        \cnodeput(4.9,1.6){F}{$\mathscr{W}_{\rho_2^*}$}
        \cnodeput(4.4,3.3){G}{$\mathscr{W}_{\rho_1^*}$}
        \psset{nodesep=0pt}
        \ncline{->}{B}{A}\lput*{:U}{$x_0x_1x_2$}
        \ncline{->}{C}{B}\lput*{:U}{$x_0x_1$}
        \ncline{->}{D}{C}\lput*{:210}{$x_0$}
        \ncline{->}{E}{D}\lput*{:180}{$x_0x_1$}
        \ncline{->}{F}{E}\lput*{:140}{$x_0$}
        \ncline{->}{G}{F}\lput*{:U}{$x_0x_1$}
        \ncline{->}{A}{G}\lput*{:U}{$x_0$}
        \ncline{->}{C}{A}\lput*{:U}{$x_1x_2x_3$}
        \ncline{->}{D}{B}\lput*{:255}{$x_3$}
        \ncline{->}{E}{C}\lput*{:200}{$x_3$}
        \ncline{->}{F}{D}\lput*{:165}{$x_3$}
        \ncline{->}{G}{E}\lput*{:110}{$x_3$}
        \ncline{->}{A}{F}\lput*{:U}{$x_3$}
        \ncline{->}{B}{G}\lput*{:U}{$x_2x_3$}
        \end{pspicture}}
      }
 \caption{Labels on arrows in the McKay quiver of type $\frac{1}{7}(1,2)$}
 \label{fig:McKayquiver}
  \end{figure}
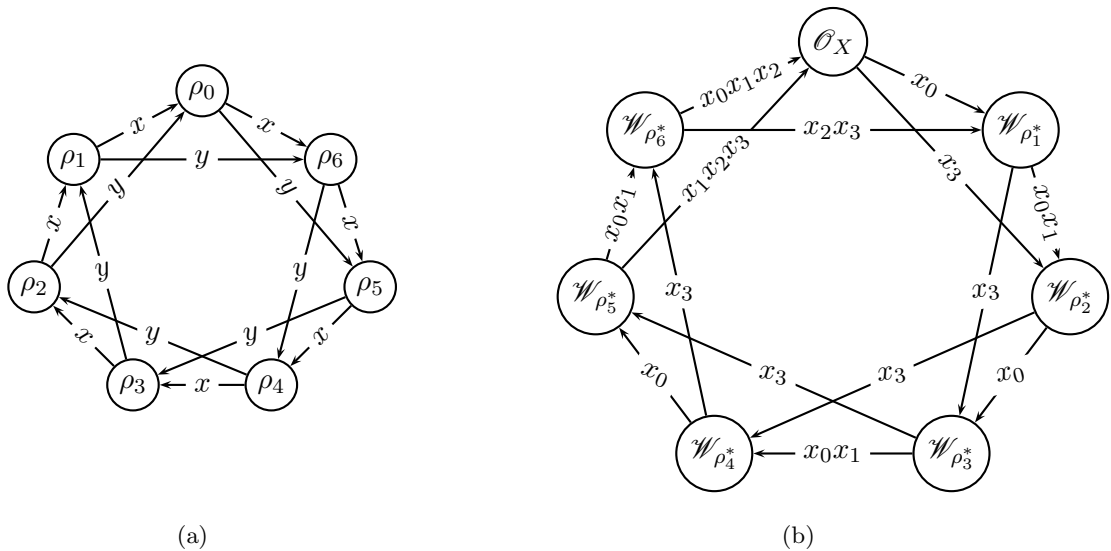
  \end{example}
  
  \subsection{The bound Special McKay quiver}
  For $\rho\in \Irr(G)$, let $(\kk[x,y]\otimes_\kk \rho^*)^G$ denote the
 $\kk[x,y]^G$-submodule of $\kk[x,y]$ generated by monomials of
 degree $\rho$.  Following  Wunram~\cite{Wunram2}, we say that $\rho^*\in \Irr(G)$ is \emph{special} if 
 $(\kk[x,y]\otimes_\kk \rho^*)^G$ has two minimal $\kk[x,y]^G$-module generators.  We define the \emph{bound Special McKay quiver} of $G\subset \GL(2,\kk)$ to be the bound quiver of sections of the sequence of line bundles $\big(\mathscr{W}_{\rho^*} : \rho^*\in \Irrsp(G)\big)$ on  $\mathcal{M}_{\vartheta^\prime}(Q^\prime,R^\prime)\cong X$ where $\Irrsp(G)$ denotes the union of the trivial representation with the special representations. 
 
 \begin{lemma}
 \label{lem:specmckay}
 The bound Special McKay quiver coincides with the bound quiver of sections $(Q,R)$ of the sequence $\underline{\mathscr{L}}$ from \eqref{eqn:sequence}. 
 \end{lemma}
 \begin{proof}
 A representation $\rho^*$ is special if and only if there is a prime exceptional divisor $D_i$ in $X$
 with toric coordinates $[x^{\alpha_i} : y^{\beta_i}]$ such that
 $\rho=\deg(x^{\alpha_i})=\deg(y^{\beta_i})$. For a special representation and for  $0\leq j\leq
 \ell$, the free $\mathscr{O}_{U_j}$-module
 $\mathscr{W}_{\rho^*}\vert_{U_j}$ of rank one is generated by
 $x^{\alpha_i}$ if $i\leq j$, and $y^{\beta_i}$ if $i > j$. 
 It follows that the tautological line bundles $\mathscr{W}_{\rho^*}$ indexed by
 special representations $\rho^*$ coincide with the line bundles $L_1,\dots,L_\ell$ satisfying $\deg(L_i\vert_{D_j}) = \delta_{ij}$ for $1\leq i,j\leq \ell$. In particular, the sequence $\big(\mathscr{W}_{\rho^*} : \rho^*\in \Irrsp(G)\big)$ coincides with the sequence $ \underline{\mathscr{L}}$ from \eqref{eqn:sequence}. 
  \end{proof}
 
  The next result describes the link between the bound McKay and Special McKay quivers.
   
 \begin{proposition}
 \label{prop:SpecialMcKay}
 The following algebras are isomorphic:
 \begin{enumerate}
 \item[\one] the quotient algebra $\kk Q/R$ for the bound Special McKay
   quiver $(Q,R)$;
 \item[\two] the endomorphism algebra
   $\End\big(\bigoplus_{0\leq i\leq \ell} L_i \big)$;
 \item[\three] the endomorphism algebra
   $\End\big(\bigoplus_{\rho^*\in \Irrsp(G)}
   \mathscr{W}_{\rho^*} \big)$;
\item[\four] the algebra $e(\kk Q^\prime/R^\prime)e$ for the bound McKay quiver $(Q^\prime,R^\prime)$ where $e=\sum_{\rho^*\in \Irrsp(G)} e_{\rho^*}$.
 \end{enumerate}
 \end{proposition}
 \begin{proof}
  Since $Q$ is the bound quiver of sections of $\underline{\mathscr{L}}$, the isomorphism
   $\one\cong\two$ follows from Lemma~\ref{lem:algebra}.   Each $L_i$ coincides with $\mathscr{W}_{\rho^*}$ for a unique $\rho^*\in
   \Irrsp(G)$, which gives the isomorphism between \two\ and \three.  Lemma~\ref{lem:boundMcKay}\two\ implies that the subalgebra $e(\kk Q^\prime/R^\prime)e$ is isomorphic to $\End\big{(}\bigoplus_{\rho^*\in \Irrsp(G)}
 \mathscr{W}_{\rho^*}\big{)}$ which gives the isomorphism $\three\cong\four$.
  \end{proof}

The assignment sending a vertex $i\in Q_0$ to the representation $\rho^*\in \Irr(G)$ for which $L_i=\mathscr{W}_{\rho^*}$ establishes a bijection $\iota\colon Q_0\to \Irrsp(G)$. Proposition~\ref{prop:SpecialMcKay} implies that for every arrow $a\in Q_1$ there is a path $p^\prime:=\iota(p)$ in $Q^\prime$ satisfying $\head(p^\prime) =\iota(\head(a))$, $\tail(p^\prime) = \iota(\tail(a))$ and $\div(p^\prime)=\div(a)$. The path $\iota(p)$ is not unique in general, though any two choices differ by
an element of the ideal of relations $R^\prime$. The situation is similar in the opposite direction:   for any path $p^\prime$ in $Q^\prime$ with $\tail(p^\prime), \head(p^\prime)\in \Irrsp(G)$, there is a path $p:=\lambda(p^\prime)$ in $Q$ satisfying $\head(p^\prime) =\iota(\head(p))$, $\tail(p^\prime) = \iota(\tail(p))$ and $\div(p^\prime)=\div(p)$; again, $\lambda(p^\prime)$ is not unique in general, though any two choices differ by
an element of the ideal $R$. In particular,  $\iota(\lambda(p^\prime))- p^\prime\in R^\prime$  and $\lambda(\iota(p))-p\in R$.

 \begin{corollary}
 For paths $p, q$ in $Q$ with the same head and tail, the following are equivalent:
 \begin{enumerate}
 \item[\one] $\div(p)=\div(q)\in \kk[x_0,\dots,x_{\ell+1}]$.
 \item[\two] $p - q \in R$;
 \item[\three] $\mon(p)=\mon(q)\in \kk[x,y]$;
 \end{enumerate}
 \end{corollary}
 \begin{proof}
 This follows from Corollary~\ref{coro:labels} and Proposition~\ref{prop:SpecialMcKay}.
  \end{proof}

 We now turn our attention to understanding the arrow set of $Q$, and to begin we cite a useful combinatorial result.  It is well known that the minimal $\kk$-algebra generators of $\kk[x,y]^G$ can be written in terms of the continued fraction expansion of $r/(r-a)$ , but it is more convenient for us to work with the expansion of $r/a$ from \eqref{eqn:ctdfrac} which encodes the nonnegative integers $c_1,\dots, c_\ell$ which in turn determine pairs $(\beta_i, \alpha_i)$ for $0\leq i\leq \ell+1$.  As Wemyss~\cite[Lemma~3.5]{Wemyss} remarks, the next combinatorial result follows from Riemenschneider's staircase~\cite{Riemenschneider1}. 
 
  \begin{lemma}
 Set $m_i:=\max \{1,c_i-1\}$ for $1\leq i\leq \ell$ and $m_{\ell+1}=1$. The inequalities 
 \begin{equation}
 \label{eqn:Riemanschneiderinequals}
 \alpha_{i+1}-m_i\alpha_i > \alpha_{i}-\alpha_{i-1}\quad \text{and}\quad m_i\beta_i-\beta_{i+1}<\beta_{i-1}-\beta_{i}
 \end{equation}
 hold for every $1\leq i\leq \ell$, and the minimal $\kk$-algebra generators of $\kk[x,y]^G$ are
 \begin{equation}
 \label{eqn:mingens}
 \{y^r\}\cup \bigcup_{0\leq i\leq \ell} \{x^{\alpha_{i+1}-t_i\alpha_i}y^{t_i\beta_i-\beta_{i+1}} : 1\leq t_i\leq m_i\}\cup \{x^r\}.
 \end{equation}  
 In particular,  listing the monomials from \eqref{eqn:mingens} with strictly decreasing exponent of $x$ is equivalent to listing the monomials with strictly increasing exponent of $y$.
 \end{lemma}
 
 Consider elements of $Q_0=\{0,1,,\dots,\ell\}$ modulo $\ell+1$, and write $\overline{i}:=i\mod \ell+1$. Recall also that $(\beta_0,\alpha_0) = (r,0)$ and $(\beta_{\ell+1},\alpha_{\ell+1}) = (0,r)$.

 \begin{proposition}
 \label{prop:xyarrows}
   For each $i\in Q_0\setminus\{0\}$ there are precisely two arrows
   with head at $i$: 
   \begin{enumerate}
   \item[\one] one arrow, denoted $a_{2i-1}$, from $i-1$ to $i$ has monomial label $\mon(a_{2i-1})= x^{\alpha_i-\alpha_{i-1}}$;
   \item[\two] the other, denoted $a_{2i+2}$, from $\overline{i+1}$ to $i$ has monomial label $\mon(a_{2i+2}) = y^{\beta_i-\beta_{i+1}}$.
   \end{enumerate}
  In particular, any arrow $a\in Q_1$ with $xy \vert \mon(a)$ has its head $\head(a)$ at vertex $0$.
 \end{proposition}
 \begin{proof}
   For $i\in Q_0\setminus\{0\}$, there are at least two arrows
   with head at $i$ since the $H^0(\mathscr{O}_X)$-algebra generators $x^{\alpha_i}, y^{\beta_i}$
   of $H^0(L_i)$ define paths in $Q$ from $0$ to $i$ with monomial a pure power of $x$ and a pure power of $y$ respectively. The path labelled $x^{\alpha_i}$ passes through vertex $i-1$, so the final arrow $a_{2i-1}$ in this path has $\mon(a_{2i-1})= x^{\alpha_i-\alpha_{i-1}}$. Similarly, the path labelled $y^{\beta_i}$ passes through $\overline{i+1}$, so the final arrow $a_{2i+2}$ in this path has label $y^{\beta_i-\beta_{i+1}}$. It remains to show that for $a\in Q_1$ with $\head(a)\in Q_0\setminus \{0\}$, the monomial $\mon(a)$ is not divisible by $xy$.   
   
   Suppose otherwise, so $\mon(a) = x^by^c$ for $b,c\geq 1$. Set $i=\head(a)$ and $j:=\tail(a)$.   The path $\iota(a)$ in the McKay quiver $Q^\prime$ from $\iota(\tail(a))$ to $\iota(\head(a))$ can be chosen to comprise $b$ arrows labelled $x$ followed by $c$ arrows labelled $y$.  If $b\geq \alpha_{j+1}-\alpha_{j}$ then this path begins with sufficiently many arrows labelled $x$ to factor via the unique path $p^\prime$ in $Q^\prime$ satisfying $\iota(a_{2j+1})=p^\prime$, but this would imply that the path $a$ in $Q$ factors via the arrow $a_{2j+1}$ which is absurd.  Also, if $c\geq \beta_{i}-\beta_{i+1}$ then $\iota(a)$ ends by traversing the unique path $q^\prime$ in $Q^\prime$ that satisfies $\iota(a_{2i+2})=q^\prime$, but this would force the path $a$ to factor via arrow $a_{2i+2}$ which is also absurd.   Similarly, the path $\iota(a)$ may also be chosen to comprise $c$ arrows labelled $y$ followed by $b$ arrows labelled $x$. If $c\geq \beta_{j-1}-\beta_{j}$ then this path in $Q^\prime$ begins with sufficiently many arrows labelled $y$ to ensure that the path $a$ in $Q$ factors via the arrow $a_{2j}$ which is absurd.  Also, if $b\geq \alpha_{i}-\alpha_{i-1}$ then the path $\iota(a)$ in $Q^\prime$ ends with sufficiently many arrows labelled $x$ to ensure that the path $a$ in $Q$ factors via the arrow $a_{2i-1}$ which is also absurd. Taken together, then, we obtain inequalities
  \begin{equation}
  \label{eqn:inequals}
   1\leq b<\min(\alpha_{j+1}-\alpha_j,\alpha_i-\alpha_{i-1})\quad\text{and}\quad 1\leq c<\min(\beta_{i}-\beta_{i+1},\beta_{j-1}-\beta_j).
  \end{equation}
 We now demonstrate that no monomial $x^by^c$ satisfying the inequalities \eqref{eqn:inequals} lies in the character space $\iota(i)\otimes\iota(j)^{-1}$. There are three cases:   
 
 \smallskip
   
   \noindent  Case 1:  $i>j$. In this case we have $\alpha_{j+1}\leq \alpha_{i}$. Multiply the generator $x^{\alpha_j}\in H^0(L_j)$ by $x^by^c$ to obtain $x^{b+\alpha_j}y^c\in H^0(L_i)$. The inequalities \eqref{eqn:inequals} give $b+\alpha_j< \alpha_{j+1}\leq\alpha_{i}$ and $c<\beta_i -\beta_{i+1} < \beta_i$. However, this implies that the section $x^{b+\alpha_j}y^c$ is divisible by neither $H^0(\mathscr{O}_X)$-algebra generator  $x^{\alpha_i}, y^{\beta_i}$ of $H^0(L_i)$ which is a contradiction.
   
   \smallskip
  
   \noindent Case 2: $i=j$. In this case, $x^by^c$ is $G$-invariant and is therefore divisible by a $G$-invariant monomial $x^{\alpha_{k+1}-t_k\alpha_k}y^{t_k\beta_k-\beta_{k+1}}$ from the list \eqref{eqn:mingens}. Then \eqref{eqn:inequals} gives $t_k\beta_k-\beta_{k+1}\leq c<\beta_{i}-\beta_{i+1}$, which forces $k>i$.   On the other hand,  inequalities \eqref{eqn:inequals} also give $\alpha_{k+1}-t_k\alpha_k\leq b< \alpha_{i}-\alpha_{i-1}$, from which we obtain $k<i$. This gives the required contradiction.
   
   \smallskip
   
   \noindent Case 3: $i<j$.  In this case we have $\beta_{j-1}\leq \beta_{i}$. Multiply the generator $y^{\beta_j}\in H^0(L_j)$ by $x^by^c$ to obtain $x^{b}y^{c+\beta_j}\in H^0(L_i)$. The inequalities \eqref{eqn:inequals} give $c+\beta_j< \beta_{j-1}\leq \beta_i$ and $b<\alpha_i -\alpha_{i-1} < \alpha_i$, but then $x^{b}y^{c+\beta_j}$ is divisible by neither $H^0(\mathscr{O}_X)$-algebra generator  $x^{\alpha_i}, y^{\beta_i}$ of $H^0(L_i)$ which is a contradiction. 
   
   \smallskip
   
\noindent   In each case we obtain a contradiction, so no such arrow exists.  This completes the proof.
\end{proof}

 \begin{remark}
 Since $x^r$ arises as the label on a cycle in $Q^\prime$ based at 0,  there exists one arrow $a$ in $Q$ from $\ell$ to 0 with $\mon(a)$ equal to a pure power of $x$. Similarly, $y^r$ labels a cycle in $Q^\prime$ based at 0, so $Q$ contains one arrow $a$  from $1$ to 0 with $\mon(a)$ equal to a pure power of $y$.
   \end{remark}

We may therefore list the arrow set as $Q_1=\{a_1,\dots, a_{\vert Q_1\vert}\}$, comprising:
\begin{itemize}
\item the \emph{$x$-arrows} $a_1,a_3,\dots, a_{2\ell+1}$, where $a_{2i-1}$ has tail at $i-1$, head at $i$ and label $x^{\alpha_{i}-\alpha_{i-1}}$;
 \item the \emph{$y$-arrows} $a_2,a_4,\dots, a_{2\ell+2}$, where $a_{2i+2}$ has tail at $\overline{i+1}$, head at $i$ and label $y^{\beta_i - \beta_{\overline{i+1}}}$;
 \item the \emph{$xy$-arrows} $a_{2\ell+3},a_{2\ell+4}\dots, a_{\vert Q_1\vert}$ each with head at 0 and $\mon(a_{2\ell+k} )$ divisible by $xy$.
  \end{itemize}
The $x$-arrows and $y$-arrows pair up naturally as $\{a_{2i+1},a_{2i+2}\}$ for $0\leq i\leq \ell$. The order in which the $xy$-arrows are listed is not important for now,  so choose any order.
 
% \begin{remark}
% Wemyss~\cite{Wemyss} proves that for each $i\in Q_0\setminus\{0\}$, the number of $xy$-arrows with tail at $i$ is given by $-(C_i^2+2)$, where $C_i^2$ is the self-intersection of the curve $C_i$. We do not 
 %\end{remark}

   \begin{example}
For the action of type $\frac{1}{21}(1,13)$, the continued fraction $21/13$ defines the pairs $(\beta_5,\alpha_5) = (0,21)$, $(\beta_4,\alpha_4) = (1,13)$, $(\beta_3,\alpha_3) = (2,5)$, $(\beta_2,\alpha_2) = (5,2)$, $(\beta_1,\alpha_1) = (13,1)$ and $(\beta_0,\alpha_0) = (21,0)$. The special representations are $\rho_{1}^*, \rho_2^*, \rho_5^*, \rho_{13}^*$ and the monomials labelling arrows in the Special McKay quiver are shown in Figure~\ref{fig:21113}.  The $\kk$-algebra generators of $\kk[x,y]^G$ from \eqref{eqn:mingens} are $y^{21}, xy^8, x^3y^3, x^8y, x^{21}$.
   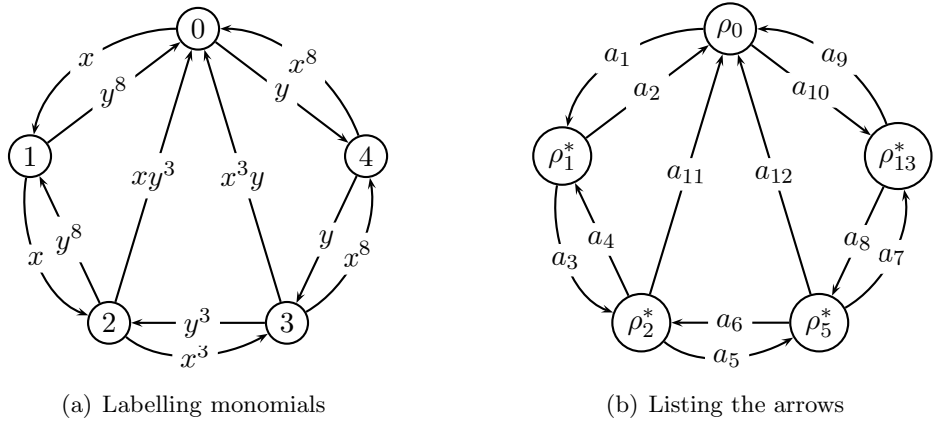
\begin{figure}[!ht]
    \centering
    \mbox{
      \subfigure[Labelling monomials]{
        \psset{unit=1.3cm}
        \begin{pspicture}(0.5,-0.5)(4.5,3.2)
        \cnodeput(2.5,3){A}{0}
        \cnodeput(0.8,1.7){B}{1} 
          \cnodeput(1.6,0){D}{2}
        \cnodeput(3.4,0){E}{3} 
          \cnodeput(4.2,1.7){G}{4}
        \psset{nodesep=0pt}
        \ncline{->}{A}{G}\lput*{:35}{$y$}
        \nccurve[angleA=0,angleB=110]{<-}{A}{G}\lput*{:35}{$x^8$}
        \nccurve[angleA=280,angleB=30]{<-}{G}{E}\lput*{:120}{$x^8$}
        \ncline{->}{G}{E}\lput*{:115}{$y$}
        \nccurve[angleA=320,angleB=210]{->}{D}{E}\lput*{:355}{$x^3$}
        \ncline{->}{E}{D}\lput*{:180}{$y^3$}
        \nccurve[angleA=260,angleB=160]{->}{B}{D}\lput*{:70}{$x$}
        \ncline{<-}{B}{D}\lput*{:65}{$y^8$}
        \nccurve[angleA=180,angleB=80]{->}{A}{B}\lput*{:145}{$x$}
        \ncline{<-}{A}{B}\lput*{:145}{$y^8$}
        \ncline{->}{D}{A}\lput*{:285}{$xy^3$}
        \ncline{->}{E}{A}\lput*{:255}{$x^3y$}
        \end{pspicture}}
      \qquad  \qquad
      \subfigure[Listing the arrows]{
        \psset{unit=1.3cm}
        \begin{pspicture}(0.5,-0.5)(4.5,3.1)
        \cnodeput(2.5,3){A}{$\rho_0$}
        \cnodeput(0.8,1.7){B}{$\rho_{1}^*$} 
          \cnodeput(1.6,0){D}{$\rho_2^*$}
        \cnodeput(3.4,0){E}{$\rho_5^*$} 
          \cnodeput(4.2,1.7){G}{$\rho_{13}^*$}
        \psset{nodesep=0pt}
        \ncline{->}{A}{G}\lput*{:35}{$a_{10}$}
        \nccurve[angleA=0,angleB=110]{<-}{A}{G}\lput*{:35}{$a_9$}
        \nccurve[angleA=280,angleB=30]{<-}{G}{E}\lput*{:120}{$a_7$}
        \ncline{->}{G}{E}\lput*{:115}{$a_8$}
        \nccurve[angleA=320,angleB=210]{->}{D}{E}\lput*{:355}{$a_5$}
        \ncline{->}{E}{D}\lput*{:180}{$a_6$}
        \nccurve[angleA=260,angleB=160]{->}{B}{D}\lput*{:70}{$a_3$}
        \ncline{<-}{B}{D}\lput*{:65}{$a_4$}
        \nccurve[angleA=180,angleB=80]{->}{A}{B}\lput*{:145}{$a_{1}$}
        \ncline{<-}{A}{B}\lput*{:145}{$a_2$}
        \ncline{->}{D}{A}\lput*{:285}{$a_{11}$}
        \ncline{->}{E}{A}\lput*{:255}{$a_{12}$}
        \end{pspicture}}     
         }
 \caption{The Special McKay quiver of type $\frac{1}{21}(1,13)$}
 \label{fig:21113}
  \end{figure}
  \end{example}

\subsection{The moduli space is irreducible}
We now prove that the fine moduli space $\mathcal{M}_\vartheta(Q,R)=\mathbb{V}(I_R)\git_\vartheta T_Q$ of $\vartheta$-stable representations of the bound Special McKay quiver is isomorphic to the image of the morphism $\varphi_{\vert\underline{\mathscr{L}}\vert}\colon X\to \vert\underline{\mathscr{L}}\vert$ from Theorem~\ref{thm:closedimmersion}. We bypass the explicit calculation of the ideal $I_R$ by introducing an auxilliary ideal $J\subset \kk[y_a : a\in Q_1]$ which suffices for our purpose. 

%As noted in Remark~\ref{rem:IRinIQ}, the toric variety $\mathbb{V}(I_Q)\git_\vartheta T_Q$ is contained in .  We now prove that $\mathbb{V}(I_Q)\git_\vartheta T_Q=\mathcal{M}_\vartheta(Q,R)$.

As a first step we exhibit a collection of elements from the ideal $R$. A cycle $p$ in $Q$ is said to be \emph{primitive} if $\mon(p)$ is one of the minimal $\kk$-algebra generators of $\kk[x,y]^G$ from \eqref{eqn:mingens}. 

\begin{lemma}
\label{lem:distinguished}
For $i\in Q_0\setminus\{0\}$, consider an arrow $a\in Q_1$ with $\tail(a)=i$. Then either:
\begin{enumerate}
\item $a=a_{2i+1}$, in which case there exists a primitive cycle $p_a$ in $Q$ based at vertex $i$  such that
$a_{2i+1}a_{2i+2} - p_a\in R$ with $a_{2i+1}, a_{2i+2}\not\in\supp(p_a)$;
\item $a=a_{2i}$, in which case there exists a primitive cycle $p_a$ in $Q$ based at vertex $i$  such that
$a_{2i} a_{2i-1} - p_a\in R$ with $a_{2i-1}, a_{2i}\not\in\supp(p_a)$;
\item $a$ is an $xy$-arrow, in which case there exists both:
\begin{enumerate}
\item[\one]  a primitive cycle $p_a$ in $Q$ based at vertex $i$ satisfying $a a_{1} a_{3}\cdots a_{2i-1} - p_a\in R$ and $a, a_{1}, a_{3}, \dots,  a_{2i-1}\not\in\supp(p_a)$; and 
\item[\two] a primitive cycle $q_a$ in $Q$ based at vertex $i$ satisfying $a a_{2\ell+2} a_{2\ell}\cdots a_{2i+2} - q_a\in R$ and $a, a_{2\ell+2}, a_{2\ell}, \dots,  a_{2i+2}\not\in\supp(q_a)$.
\end{enumerate}
\end{enumerate}
\end{lemma}
\begin{proof}
 For case (1),  set $a=a_{2i+1}$. Since $a$ is an $x$-arrow and $a_{2i+2}$ is a $y$-arrow, Proposition~\ref{prop:SpecialMcKay} implies that the cycle $a a_{2i+2}$ in $Q$ based at $i$ arises from a cycle $p^\prime:= \iota(a a_{2i+2})$ in $Q^\prime$ based at $\iota(i)\in Q_0^\prime$ where $\mon(p^\prime)=x^my^n$ with $m,n\geq 1$. Since $a$ precedes $a_{2i+2}$, the first arrow $a_1^\prime\in Q_1^\prime$ traversed by $p^\prime$ satisfies $\mon(a_1^\prime) = x$. Define $p^\prime_a$ in $Q^\prime$ to be the unique cycle based at $\iota(i)\in Q_0^\prime$ with $\mon(p_a^\prime) = x^my^n$ that first traverses $n$ arrows labelled $y$ and then traverses $m$ arrows labelled $x$. The unique cycle $p_a$ in $Q$ satisfying $p_a^\prime=\iota(p_a)$ first traverses an arrow $a_f$ with $\tail(a_f)=i$ and $y\vert \mon(a_f)$ and ends by traversing an arrow $a_l$ with $\head(a_l)=i$ and $x\vert \mon(a_l)$. Proposition~\ref{prop:xyarrows} implies that $a_l=a_{2i-1}$. If $a_f$ is a $y$-arrow then $a_f=a_{2i}$ and $p_a = a_{2i}a_{2i-1}$ is the only such primitive cycle; otherwise, $a_f$ is an $xy$-arrow, in which case $p_a = a_f a_{1} a_{3}\cdots a_{2i-1}$ is the only such primitive cycle. In either case $a_{2i+1}, a_{2i+2}\not\in\supp(p_a)$, so case (1) is complete. Case (2) is the same as case (1) with the roles of $x$ and $y$ switched.  For case (3),  let $a$ be an $xy$-arrow, so $\tail(a)=i$, $\head(a)=0$ and $\mon(a)=x^my^n$ for $m, n\geq 1$. The cycle $p^\prime:= \iota(a a_{1} a_{3}\cdots a_{2i-1})$ in $Q^\prime$ based at $\iota(i)\in Q_0^\prime$ ends with at least $\alpha_i$ arrows labelled $x$ and satisfies $\mon(p^\prime)=x^{m+\alpha_i}y^n$. Define $p^\prime_a$ in $Q^\prime$ to be the unique cycle based at $\iota(i)\in Q_0^\prime$ with $\mon(p_a^\prime) = x^{m+\alpha_i}y^n$ that first traverses $m+\alpha_i$ arrows labelled $x$ and then traverses $n$ arrows labelled $y$.  The unique cycle $p_a$ in $Q$ satisfying $p_a^\prime=\iota(p_a)$ first traverses an arrow $a_f$ with $\tail(a_f)=i$ and $x\vert \mon(a_f)$, and ends by traversing an arrow $a_l$ with $\head(a_l)=i$ and $y\vert \mon(a_l)$. Proposition~\ref{prop:xyarrows} gives $a_l=a_{2i+2}$.  If $a_f$ is an $x$-arrow then $a_f=a_{2i+1}$ and $p_a = a_{2i+1}a_{2i+2}$ is the only such primitive cycle; otherwise, $a_f$ is an $xy$-arrow, in which case $p_a = a_f a_{2\ell+2} a_{2\ell}\cdots a_{2i+2}$ is the only such primitive cycle.  Note that $a_f\neq a$ because $x^{\alpha_i}=\mon(a_{1} a_{3}\cdots a_{2i-1})\neq \mon(a_{2\ell+2} a_{2\ell}\cdots a_{2i+2})=y^{\beta_i}$.  It follows that in either case $a, a_{1}, a_{3}, \dots,  a_{2i-1}\not\in\supp(p_a)$, so case (3.i) is complete. Case (3.ii) is similar.
 \end{proof}

Let $\Lambda$ denote the set of all path differences constructed in Lemma~\ref{lem:distinguished} as the vertex $i$ ranges over $Q_0\setminus\{0\}$ and the arrow $a$ ranges over $\{a\in Q_1 : \tail(a)=i\}$. Define the ideal
 \[
 J:=\big( y_p-y_q\in \kk[y_a : a\in Q_1] : p-q\in\Lambda\big).
 \]
 Remark~\ref{rem:IRinIQ} and Lemma ~\ref{lem:distinguished} imply that $J\subseteq I_R\subseteq I_Q$, and hence $\mathbb{V}(I_Q)\subseteq \mathbb{V}(I_R) \subseteq \mathbb{V}(J)$. 
 
 \begin{remark}
 We make no attempt to prove that $J=I_R$ since this requires an understanding of $R$ itself.  The explicit calculation of $R$ by Wemyss~\cite[pages 3-14]{Wemyss} does in fact imply that $J=I_R$, but we emphasise that we do not make use of this observation. Rather, our approach is to show that while $\mathbb{V}(J)\neq \mathbb{V}(I_Q)$, we nevertheless get equality after removing the locus cut out by the irrelevant ideal, that is,  $\mathbb{V}(J)\setminus \mathbb{V}(B_Q) =\mathbb{V}(I_R)\setminus \mathbb{V}(B_Q) =  \mathbb{V}(I_Q)\setminus \mathbb{V}(B_Q)$. We illustrate this with a pair of examples.
   \end{remark}
 
\begin{example}
For the action of type $\frac{1}{7}(1,2)$, the quiver $Q$ from Figure~\ref{fig:712}(b) is obtained from the McKay quiver $Q^\prime$ from Figure~\ref{fig:McKayquiver}(b) by keeping only paths with head and tail at the vertices indexed by $\mathscr{O}_X=\mathscr{W}_{\rho_0^*}$, $L_1 = \mathscr{W}_{\rho_1^*}$ and $L_2 = \mathscr{W}_{\rho_2^*}$.   To obtain the monomials $\mon(a)$ for $a\in Q_1$, consider the labels in Figure~\ref{fig:712}(b) and replace $x_0$ by $x$, $x_3$ by $y$ and set $x_1=x_2=1$. Lemma~\ref{lem:distinguished} constructs the ideal $J$ as follows: arrows $a_2$ and $a_3$ with tail at $1$ each determine $y_1y_2-y_3y_4\in \Lambda$; the $x$-arrow $a_5$ determines $y_5y_6-y_1y_3y_7$; the $y$-arrow $a_4$ determines $y_3y_4-y_6y_8$;  the $xy$-arrow $a_7$ determines both $y_1y_3y_8-y_6y_7$ and $y_6y_8-y_3y_4$; and the $xy$-arrow $a_8$ determines both $y_6y_7-y_1y_3y_8$ and $y_1y_3y_7-y_5y_6$. Taken together we obtain 
   \begin{eqnarray*}
   J & = & \big(y_1y_2-y_3y_4, \: y_5y_6-y_1y_3y_7,\: y_3y_4-y_6y_8,\: y_1y_3y_8-y_6y_7\big), \\
      & = & I_Q \; \cap (y_1, y_4, y_6)\:\cap\: (y_1, y_3, y_6)\:\cap\; (y_2, y_3, y_6).
   \end{eqnarray*}
 Saturating by the irrelevant ideal $B_Q = (y_1y_3, y_1y_6, y_4y_6)$ gives $J : B_Q=I_Q : B_Q$, and hence $\mathbb{V}(J)\setminus \mathbb{V}(B_Q) = \mathbb{V}(I_Q)\setminus \mathbb{V}(B_Q)$.
\end{example}
 
  \begin{example}
The Special McKay quiver for the action of type $\frac{1}{21}(1,13)$ is shown in Figure~\ref{fig:21113}.   Applying Lemma~\ref{lem:distinguished} repeatedly gives
 \[
 J=\left(\begin{array}{c}  y_1y_2-y_3y_4,\: y_7y_8-y_9y_{10}, \: y_8y_{10}y_{12}-y_5y_6, \\ y_1y_3y_{11}-y_5y_6, \:y_1y_3y_5y_{12}-y_9y_{10},\: y_6y_8y_{10}y_{11}-y_3y_4\end{array}\right).
 \]
This ideal has 25 primary components, one of which is the toric ideal $I_Q$ and 24 of which cut out linear varieties.  Saturating by the irrelevent ideal 
 \[
 B_Q = (y_1y_3y_5y_7,\: y_1y_3y_5y_{10}, \:y_1y_3y_8y_{10},\: y_1y_6y_8y_{10},\: y_4y_6y_8y_{10})
 \]
 gives $J : B_Q = I_Q : B_Q$, from which we obtain $\mathbb{V}(J)\setminus \mathbb{V}(B_Q) = \mathbb{V}(I_Q)\setminus \mathbb{V}(B_Q)$.
 \end{example}

 \begin{theorem}
 The fine moduli space $\mathcal{M}_\vartheta(Q,R)$ of $\vartheta$-stable representations of the bound Special McKay quiver is isomorphic to the toric variety $\mathbb{V}(I_Q)\git_\vartheta T_Q$.
 \end{theorem}
 \begin{proof}
 Proposition~\ref{prop:generic} implies that $\mathbb{A}^{Q_1}_\kk\setminus \mathbb{V}(B_Q)$ is covered by charts $\Spec\big( \kk[y_a : a\in Q_1]_{y_\mathcal{T}}\big)$, where $\mathcal{T}\subseteq Q_1$ is a spanning tree with root at $0\in Q_0$ and $y_{\mathcal{T}}:=\prod_{a\in \supp(\mathcal{T})} y_a$.  The inclusion $J\subseteq I_Q$ implies that $\mathbb{V}(I_Q)\setminus \mathbb{V}(B_Q)$ is a closed subscheme of $\mathbb{V}(J)\setminus \mathbb{V}(B_Q)$, and restricting to the open cover gives $\Spec\big( (\kk[y_a : a\in Q_1]/I_Q)_{y_\mathcal{T}}\big) \subseteq \Spec\big( (\kk[y_a : a\in Q_1]/J)_{y_\mathcal{T}}\big)$ for each $\mathcal{T}$. The variety $\mathbb{V}(I_Q)\git_\vartheta T_Q$ is the geometric quotient of $\mathbb{V}(I_Q)\setminus \mathbb{V}(B_Q)$ by the action of $T_Q$, and since each toric chart of $\mathbb{V}(I_Q)\git_\vartheta T_Q\cong X$ is isomorphic to $\mathbb{A}^2_\kk$, it follows that $\Spec\big( (\kk[y_a : a\in Q_1]/I_Q)_{y_\mathcal{T}}\big)$ is isomorphic to $T_Q\times \mathbb{A}^2_\kk$ for every spanning tree $\mathcal{T}$. We claim that 
 $\Spec\big( (\kk[y_a : a\in Q_1]/J)_{y_\mathcal{T}}\big)$ is also isomorphic to $T_Q\times \mathbb{A}^2_\kk$ for every $\mathcal{T}$, giving $\mathbb{V}(I_Q)\setminus \mathbb{V}(B_Q) = \mathbb{V}(J)\setminus \mathbb{V}(B_Q)$. The result follows from the claim, because the inclusions $J\subseteq I_R\subseteq I_Q$ force $\mathbb{V}(I_R)\setminus \mathbb{V}(B_Q) = \mathbb{V}(I_Q)\setminus \mathbb{V}(B_Q)$, and the required identification $\mathcal{M}_\vartheta(Q,R):= \mathbb{V}(I_R)\git_\vartheta T_Q= \mathbb{V}(I_Q)\git_\vartheta T_Q$ follows from Proposition~\ref{prop:generic}.

To prove the claim we first list all spanning trees $\mathcal{T}$ in $Q$ with root at $0\in Q_0$.  Proposition~\ref{prop:xyarrows} implies that there are $\ell+1$ such trees, namely,  $\mathcal{T}_0,\dots, \mathcal{T}_{\ell}$ where $\supp(\mathcal{T}_0) = \{a_4,a_6,\dots,a_{2\ell+2}\}$, $\supp(\mathcal{T}_j) = \{a_{1}, a_{3}, \dots, a_{2j-1}\}\cup\{a_{2j+4}, a_{2j+6},\dots, a_{2\ell+2},\}$ for every $0 < j <\ell$, and $\supp(\mathcal{T}_\ell) = \{a_1,a_3,\dots,a_{2\ell-1}\}$. The result follows once we prove that 
 \begin{equation}
 \label{eqn:localisation}
 \Big(\frac{\kk[y_a : a\in Q_1]}{J}\Big)_{y_{\mathcal{T}_j}}\cong \kk[y_a^{\pm 1} : a\in \supp(\mathcal{T}_j)]\otimes_\kk \kk[y_{2j+1}, y_{2j+2}]\quad\text{for }0\leq j\leq \ell.
 \end{equation}
Fix $0\leq j\leq \ell$. We use the ideal $J$ to eliminate certain variables $y_a$ in two stages (if $j=0$ ignore stage one, and if $j=\ell$ ignore stage two): 
\begin{enumerate}
 \item \emph{eliminate each $y_a$ with $1\leq \tail(a)\leq j$ and $y\vert\mon(a)$}.  Fix $i\in Q_0\setminus\{0\}$ with $i\leq j$ and consider the $y$-arrow $a_{2i}$. Lemma~\ref{lem:distinguished}(2) gives $y_{2i} = y_{p_a} y_{2i-1}^{-1}$ in $\big(\kk[y_a : a\in Q_1]/J\big)_{y_T}$ where $y_{p_a}$ is either $y_{2i+1}y_{2i+2}$ or $y_{a_f}y_{2\ell+2}\dots y_{2i+2}$ for some $xy$-arrow $a_f$ with $\tail(a_f)=i$. In either case,  $y_{2i}$ is expressed in terms of  variables with a higher index, and we eliminate it.  As for the $xy$-arrows with tail at $i$, choose one such with $\mon(a)$ divisible by the lowest power of $x$. Lemma~\ref{lem:distinguished}(3.i) gives $y_a = y_{p_a}\cdot(y_{1}y_{3}\cdots y_{2i-1})^{-1}$ in $\big(\kk[y_a : a\in Q_1]/J\big)_{y_T}$, where $y_{p_a}$ is either $y_{2i+1}y_{2i+2}$ or $y_{a_f}y_{2\ell+2}\dots y_{2i+2}$ for some $xy$-arrow $a_f$ with tail at $i$ and $\mon(a_f)$ is divisible by a strictly higher power of $x$ than is $\mon(a)$. In either case,  $y_a$ is expressed in variables indexed by arrows from $\{a_1,\dots,a_{2i+ 1}\}\cup\{a_{2i+2},\dots,a_{2\ell+2}\}$ and $xy$-arrows $b$ with tail at $i$ and $\mon(b)$ divisible by a strictly higher power of $x$ than was $\mon(a)$. Repeat, until each $y_a$ indexed by an $xy$-arrow $a$ with tail at $i$ is expressed in variables $\{y_1,\dots,y_{2i+1}\}\cup\{y_{2i+2},\dots,y_{2\ell+2}\}$.  Carry out this procedure, in order, at the vertices on the list $(1,2,\dots,j)$. The result is that each $y_a$ with $1\leq \tail(a)\leq j$ and $y\vert\mon(a)$ is written in terms of $\{y_{1},y_3,\dots,y_{2j+ 1}\}\cup\{y_{2j+2},\dots,y_{2\ell+2}\}$.

\item \emph{eliminate each $y_a$ with $j+1\leq \tail(a)\leq \ell$ and $x\vert\mon(a)$}.  Fix $i\geq j+1$ and consider the $x$-arrow $a_{2i+1}$. Lemma~\ref{lem:distinguished}(1) gives $y_{2i+1} = y_{p_a} y_{2i+2}^{-1}$ in $\big(\kk[y_a : a\in Q_1]/J\big)_{y_T}$ where $y_{p_a}$ is either $y_{2i-1}y_{2i}$ or $y_{a_f}y_{1}\dots y_{2i-1}$ for some $xy$-arrow $a_f$ with tail at $i$. In either case,  $y_{2i+1}$ is expressed in terms of $x$-arrows $\{a_1, a_{3}\cdots a_{2i-1}\}$, the $y$-arrows $\{a_{2i}, a_{2i+2}\}$ and an $xy$-arrow with tail at $i$, so we eliminate it.  As for the $xy$-arrows with tail at $i$, choose one such with $\mon(a)$ divisible by the lowest power of $y$. Lemma~\ref{lem:distinguished}(3.ii) gives $y_a = y_{q_a}\cdot(y_{2\ell+2}y_{2\ell}\cdots y_{2i+2})^{-1}$ in $\big(\kk[y_a : a\in Q_1]/J\big)_{y_T}$, where $y_{q_a}$ is either $y_{2i-1}y_{2i}$ or $y_{a_f}y_{1}\dots y_{2i-1}$ for some $xy$-arrow $a_f$ with tail at $i$ and $\mon(a_f)$ is divisible by a strictly higher power of $y$ than is $\mon(a)$.  In either case,  $y_a$ is expressed in variables indexed by arrows from $\{a_1,\dots,a_{2i-1}\}\cup\{a_{2i},\dots,a_{2\ell+2}\}$ and $xy$-arrows $b$ with tail at $i$ and $\mon(b)$ divisible by a higher power of $y$ than was $\mon(a)$.  Repeat, until each $y_a$ indexed by an $xy$-arrow $a$ with tail at $i$ is expressed in variables $\{y_1,\dots,y_{2i-1}\}\cup\{y_{2i},\dots,y_{2\ell+2}\}$.  Carry out this procedure, in order, at the vertices on the list $(\ell, \ell-1,\dots,j+1)$, so each $y_a$ with $j+1\leq \tail(a)\leq \ell$ and $x\vert\mon(a)$ is expressed in $\{y_{1},y_3,\dots,y_{2j+ 1}\}\cup\{y_{2j+2},\dots,y_{2\ell+2}\}$. 
\end{enumerate}
 
Applying Stages 1 and 2 gives each variable $y_a$ from $\big(\kk[y_a : a\in Q_1]/J\big)_{y_T}$ in terms of the variables $\{y_a : a\in \supp(\mathcal{T}_j)\}\cup\{y_{2j+1},y_{2j+2}\}$. It follows that the $\kk$-algebra homomorphism
\[
\kk[y_a^{\pm 1} : a\in \supp(\mathcal{T}_j)]\otimes_\kk \kk[y_{2j+1}, y_{2j+2}] \longrightarrow \big(\kk[y_a : a\in Q_1]/J\big)_{y_T}
\]
 sending $y_a$ to $[y_a\mod J]$ is an epimorphism. If the kernel is nonzero then the affine chart  $\Spec\big( (\kk[y_a : a\in Q_1]/J)_{y_\mathcal{T}}\big)$ is a proper subscheme of $T_Q\times \mathbb{A}^2_\kk\cong \Spec\big(\kk[y_a^{\pm 1} : a\in \supp(\mathcal{T}_j)]\otimes_\kk \kk[y_{2j+1}, y_{2j+2}]\big)$, contradicting the inclusion $\mathbb{V}(I_Q)\setminus \mathbb{V}(B_Q)\subseteq \mathbb{V}(J)\setminus \mathbb{V}(B_Q)$.
  \end{proof}
  
 \section{The Special McKay correspondence following Van den Bergh}
 
 To conclude we consider the derived category.  Motivated by Bridgeland's work on perverse coherent sheaves, Van den
 Bergh~\cite{VDB} considered the subcategory $\mathcal{B}:=\Per$ of the bounded derived category of coherent sheaves on $X$ consisting
 of those objects $E$ whose cohomology sheaves $\mathcal{H}^i(E)$ are
 nonzero only in degrees $-1$ and $0$, such that
 $f_*(\mathcal{H}^{-1}(E))=0$, and such that the map
 $f^*f_*\mathcal{H}^0(E)\to \mathcal{H}^0(E)$ is surjective. A
 \emph{projective generator} for $\mathcal{B}$ is defined to be a
 projective object $P$ in $\mathcal{B}$ for which
 $\Hom_{\mathcal{B}}(P,E)=0$ implies $E=0$.

 \begin{lemma}
 \label{lem:progen}
 The vector bundle $\bigoplus_{i\in Q_0} L_i$ on $X$ is a projective generator
 for $\mathcal{B}$.
 \end{lemma}
 \begin{proof}
   The minimal resolution satisfies $\textbf{R}f_*(\mathscr{O}_X) =
   \mathscr{O}_{\mathbb{A}^2_\kk/G}$.  Since $\kk$ is algebraically
   closed of characteristic zero, and since each fibre of $f$ has
   dimension at most one, the results from Van den
   Bergh~\cite[Section~3]{VDB} hold. In particular, since the ample
   line bundle $L=\bigotimes_{i=1}^\ell L_i$ satisfies
   $L=\det(\bigoplus_{i=1}^\ell L_i)$, the Artin--Verdier construction
   gives a short exact sequence
 \[
  \begin{CD}   
    0 @>>> \mathscr{O}_X^{\oplus(\ell-1)} @>>> \textstyle{\bigoplus_{i=1}^\ell L_i} @>>> L @>>> 0
  \end{CD}
 \]
 by \cite[Lemma~3.5.1]{VDB}.  Since $L$ is globally
 generated,  the proof of \cite[Proposition~3.2.5]{VDB} shows that
 the bundle $\mathscr{O}_X\oplus \bigoplus_{i=1}^\ell L_i$ is a
 projective generator for $\mathcal{B}$.
 \end{proof}

 For $A=\End\bigl( \bigoplus_{i\in Q_0} L_i
 \bigr)$, write $\modA$ for the category of
 finitely generated right $A$-modules.  

 \begin{proof}[Proof of Theorem~\ref{thm:mainderived}]
 Since $\bigoplus_{i\in Q_0} L_i$ on $X$ is a projective generator for
 $\mathcal{B}$, the construction \cite[Corollary~3.2.8]{VDB}
 gives an equivalence between $\mathcal{B}$ and $\modA$ that extends
 to an equivalence 
  \[
\textbf{R}\Hom\bigl(\textstyle{\bigoplus_{i\in Q_0}} L_i, -\bigr)
 \colon D^b\big(\!\Coh(X)\big)\longrightarrow D^b\big(\!\modA\big)
 \]
 of derived categories. This completes the proof of Theorem~\ref{thm:mainderived}.
 \end{proof}

% \bibliography{specialMcKay}

\begin{thebibliography}{10}

\bibitem{BKR}
T.~Bridgeland, A.~King, and M.~Reid.
\newblock The {M}c{K}ay correspondence as an equivalence of derived categories.
\newblock {\em J. Amer. Math. Soc.}, 14(3):535--554 (electronic), 2001.

\bibitem{Cox}
D.~Cox.
\newblock The homogeneous coordinate ring of a toric variety.
\newblock {\em J. Algebraic Geom.}, 4(1):17--50, 1995.

\bibitem{CMT1}
A.~Craw, D.~Maclagan, and R.R. Thomas.
\newblock Moduli of {M}c{K}ay quiver representations {I}: the coherent
  component.
\newblock {\em Proc. London Math. Soc.}, \textbf{95}(1):179--198, 2007.

\bibitem{CMT2}
A.~Craw, D.~Maclagan, and R.R. Thomas.
\newblock Moduli of {M}c{K}ay quiver representations. {II}. {G}r\"obner basis
  techniques.
\newblock {\em J. Algebra}, 316(2):514--535, 2007.

\bibitem{CrawSmith}
A.~Craw and G.~G. Smith.
\newblock Projective toric varieties as fine moduli spaces of quiver
  representations.
\newblock {\em Amer. J. Math.}, 130(6):1509--1534, 2008.

\bibitem{M2}
D.~Grayson and M.~Stillman.
\newblock Macaulay 2, a software system for research in algebraic geometry.
\newblock Available from {\tt http://www.math.uiuc.edu/Macaulay2/}.

\bibitem{Ishii}
A.~Ishii.
\newblock On the {M}c{K}ay correspondence for a finite small subgroup of {${\rm
  GL}(2,\mathbb C)$}.
\newblock {\em J. Reine Angew. Math.}, 549:221--233, 2002.

\bibitem{ItoNakamura}
Y.~Ito and I.~Nakamura.
\newblock Hilbert schemes and simple singularities.
\newblock In {\em New trends in algebraic geometry (Warwick, 1996)}, volume 264
  of {\em London Math. Soc. Lecture Note Ser.}, pages 151--233. Cambridge Univ.
  Press, Cambridge, 1999.

\bibitem{KapranovVasserot}
M.~Kapranov and E.~Vasserot.
\newblock Kleinian singularities, derived categories and {H}all algebras.
\newblock {\em Math.~Ann.}, 316:565--576, 2000.

\bibitem{Kidoh}
R~Kidoh.
\newblock Hilbert schemes and cyclic quotient surface singularities.
\newblock {\em Hokkaido Math. J.}, \textbf{30}(1):91--103, (2001).

\bibitem{King}
A.~King.
\newblock Moduli of representations of finite-dimensional algebras.
\newblock {\em Quart.~J.~Math.~Oxford Ser.~(2)}, \textbf{45}(180):515--530,
  1994.

\bibitem{Riemenschneider1}
O.~Riemenschneider.
\newblock Deformationen von {Q}uotientensingularit\"aten (nach zyklischen
  {G}ruppen).
\newblock {\em Math. Ann.}, 209:211--248, 1974.

\bibitem{VDB}
M.~Van~den Bergh.
\newblock Three-dimensional flops and noncommutative rings.
\newblock {\em Duke Math. J.}, 122(3):423--455, 2004.

\bibitem{Wemyss}
M.~Wemyss.
\newblock Reconstruction algebras of type {$A$}, 2007.
\newblock \texttt{arXiv:0704.3693}, to appear \emph{Trans.\ Amer.\ Math.\ Soc.}

\bibitem{Wunram2}
J.~Wunram.
\newblock Reflexive modules on cyclic quotient surface singularities.
\newblock In {\em Singularities, representation of algebras, and vector bundles
  (Lambrecht, 1985)}, volume 1273 of {\em Lecture Notes in Mathematics}, pages
  221--231. Springer, Berlin, 1987.

\bibitem{Wunram}
J.~Wunram.
\newblock Reflexive modules on quotient surface singularities.
\newblock {\em Math. Ann.}, \textbf{279}(4):583--598, (1988).

\end{thebibliography}

 \end{document}